\documentclass[final]{siamltex} 
\usepackage{graphicx}
\usepackage{graphicx,amssymb,amsmath}
\usepackage{epsfig}
\usepackage{subfigure}
\usepackage{cite}
\usepackage{booktabs}
\usepackage{array}
\usepackage{color}
\usepackage{enumerate}
\usepackage{algorithm}
\usepackage{algorithmic}
\usepackage{float}
\usepackage{multirow}
\usepackage{rotating}
\usepackage[misc]{ifsym}
\usepackage{bm}
\usepackage{mathrsfs}
\usepackage[right=1.1in]{geometry}
\usepackage{mathrsfs}
\usepackage[misc]{ifsym} 
\usepackage{caption} 
\newtheorem{exmp}{Example}

\begin{document}

\title{Structured FISTA for Image Restoration}

\author{
Zixuan Chen\thanks{School of Mathematical Sciences, Dalian University of Technology, Dalian 116025, Liaoning, China({\tt chenzixuan@mail.dlut.edu.cn}, {\tt yubo@dlut.edu.cn}).}
\and
James G. Nagy\thanks{Department of Mathematics, Emory University, Atlanta, GA 30322, U.S.A
({\tt jnagy@emory.edu}, {\tt yuanzhe.xi@emory.edu}). This research was supported in part by the U.S. National Science Foundation under grant DMS-1819042}
\and
Yuanzhe Xi\footnotemark[2]
\and Bo Yu\footnotemark[1]
}

\maketitle
\begin{abstract}
In this paper, we propose an efficient numerical scheme for solving some large scale ill-posed linear inverse problems arising from image restoration.  In order to accelerate the computation, two different hidden structures are exploited. First, the coefficient matrix is approximated as the sum of a small number of Kronecker products. This procedure not only introduces one more level of parallelism into the computation but also enables the usage of computationally intensive matrix-matrix multiplications in the subsequent optimization procedure. We then derive the corresponding Tikhonov regularized minimization model and extend the fast iterative shrinkage-thresholding algorithm (FISTA) to solve the resulting optimization problem. Since the matrices appearing in the Kronecker product approximation are all structured matrices (Toeplitz, Hankel, etc.), we can further exploit their fast matrix-vector multiplication algorithms at each iteration. The proposed algorithm is thus called \emph{structured fast iterative shrinkage-thresholding algorithm} (sFISTA). In particular, we show that the approximation error introduced by sFISTA is well under control and sFISTA can reach the same image restoration accuracy level as FISTA. Finally, both the theoretical complexity analysis and some numerical results are provided to demonstrate the efficiency of sFISTA.
\end{abstract}

% REQUIRED
\begin{keywords}
  linear inverse problem, image restoration, Kronecker product approximation, structured FISTA
\end{keywords}

% REQUIRED
%\begin{AMS}
%  68Q25, 68R10, 68U05
%\end{AMS}

\section{Introduction}\label{sec: Introduction}
Image restoration problems have a wide range of important applications, such as digital camera and video, microscopy, meidcal imaging, etc.. Image restoration is the process of reconstructing an image of an unknown scene from an observed image, where the distortion can arise from many sources, such as motion blurs, out of focus lens, or atmospheric turbulence. Suppose there is an exact image of being all black except for a single bright pixel. If we take a picture of this image, then the distortion operation will cause the single bright pixel to be spread over its neighboring pixels. This single bright pixel is called a point source, and the function that describes the distortion and the resulting image of the point source is called the point spread function (PSF) \cite{Hansen2006D}. Mathematically, the distortion can be represented by a PSF. If the PSF is the same regardless of the location of the point source, it is called spatially invariant. Throughout this paper, we assume the PSF under consideration is always spatially invariant.

A spatially invariant image restoration problem can be modeled by a linear inverse problem of the following form
\begin{equation}\label{linear system}
b=Ax+e,
\end{equation}
where $A\in{\mathbb{R}}^{N\times N}$ is a blurring matrix constructed from the PSF, $e\in{\mathbb{R}}^{N}$ is a vector representing additive noise, $b\in{\mathbb{R}}^{N}$ represents the distorted image and $x\in{\mathbb{R}}^{N}$ denotes the unknown true image to be estimated. The matrix $A$ is usually very ill-conditioned in these image restoration problems. 

A classical way to solve \eqref{linear system} is by the least squares (LS) approach\cite{Bjorck1996N}, whose solution takes the following form
\begin{equation*}\label{LS model}
\hat{x}_{\mathrm{LS}}=\arg\min_{x}\ \frac{1}{2}\|Ax-b\|_2^2.       
\end{equation*}
However, when $A$ is ill-conditioned, the LS solution usually has a huge norm and is thus meaningless\cite{Hansen2006D}. In order to compute a decent approximation to $x$, it is necessary to employ some form of regularization. The basic idea of regularization is to replace the original ill-conditioned problem with a ``nearby" well-conditioned problem whose solution is close to the orignal solution. Tikhonov regularization \cite{Tikhonov1977} is one of the most popular regularization techniques, where a quadratic penalty is added to the object function
\begin{equation*}\label{TIK model}
\hat{x}_{\mathrm{TIK}}=\arg\min_{x}\ \frac{1}{2}\|Ax-b\|_2^2+\frac{{\lambda}^{2}}{2}\|Rx\|_2^2.        
\end{equation*}
%\JN{Change $L$ for regularization matrix to $R$ so it won't get confused later with the Lagrange multiplier.}
The second term in the above equation is a regularization term, which controls the norm (or seminorm) of the solution. The reguarization parameter $\lambda\!>\!0$ controls ``smoothness" of the regularized solution. Typical choices of $R$ include an identity matrix and a matrix approximating the first or second order derivative operator \cite{Golub1999T, Hansen1993T, Hansen1997R}.

In this paper, we choose $R$ as an identity matrix and consider the following minimization model
\begin{equation}\label{vector minimization model}
\min_{x}\ {\Phi(x)}=\frac{1}{2}\|Ax-b\|_2^2+\frac{{\lambda}^{2}}{2}\|x\|_2^2%\tag{{$\mathcal{P}$}}.        
\end{equation}
In many applications, such as image restoration, it may also be important to include convex constraints (e.g., $x \geq 0$) on the solution.

{Numerous algorithms proposed in the literature can be used to solve \eqref{vector minimization model} with convex constraints. One of them is the interior point method \cite{Ben2001, nocedal2006numerical}. However, image restoration problems often involve dense matrix data, which will hamper the effectiveness of the interior point method. Another popular class of methods for solving \eqref{vector minimization model} are gradient-based algorithms \cite{Beck2009Fast, Beck2009Gradient, Zhang2018A}. Although these algorithms are relatively inexpensive at each iteration, they often suffer from slow convergence. One recent development is the fast iterative shrinkage-thresholding algorithm (FISTA) \cite{Beck2009A}, which was proposed to solve nonsmooth convex optimization problems. FISTA preserves the computational simplicity and has a fast global convergence rate. Thus, FISTA becomes quite attractive for solving large-scale problems. Although problem \eqref{vector minimization model} does not involve any nonsmooth term, 
incorporating convex constraints is important in image deblurring applications. Moreover, in some situations $l_1$-based regularization has to be exploited to enforce sparsity in the solution.
We plan to apply the proposed method to solve this class of nonsmooth optimization problems in the future. In this paper, we will first fully take advantage of the hidden structures of the blurring matrix $A$ and  improve the efficiency of the FISTA framework for solving the smooth optimization problem \eqref{vector minimization model}.
}

Since the blurring model is essentially a convolution, the first structure to be exploited is the Kronecker product structure. Assume $K\in{\mathbb{R}}^{n\times n}$ and $H\in{\mathbb{R}}^{m\times m}$, the Kronecker product of these two matrices is defined as
\begin{equation}\label{Kronecker product}
K\otimes H=\left[
\begin{array}{ccc}
k_{11}H&\cdots&k_{1n}H\\
\vdots&&\vdots\\
k_{n1}H&\cdots&k_{nn}H
\end{array}\right].        
\end{equation}

For the blurring operator $A$ in \eqref{linear system}, it has been shown that $A$ can be approximated by a matrix $A_s$ as follows \cite{Kamm2000, Nagy2004}
\begin{equation}\label{Kronecker structure}
A\approx A_s=\sum\limits_{i=1}^{s}\ K_i\otimes H_i,
\end{equation}
where $K_i\in{\mathbb{R}}^{n\times n}$, $H_i\in{\mathbb{R}}^{m\times m}$ with $N=mn$.  
{The error between the blurring matrix and the Kronecker product approximation can be easily controlled. In addition, these $K_i$ and $H_i$ are not general dense matrices but structured matrices (Toeplitz, Hankel, etc. \cite{Xi2014A, Xi2014S, Xia2012A}). We will give more details on the error between $A$ and $A_s$ and the structures of $K_i$ and $H_i$ in Section \ref{sec: Kroneker decomposition}.} 

Consequently, the solution of \eqref{linear system} can be approximated by the following problem
\begin{equation}\label{linear system approx}
b=A_sx_s+e
\end{equation} 
and equivalently the solution of \eqref{vector minimization model} can be approximated by solving the optimization problem
\begin{equation}\label{equivalent vector minimization model}
\min_{x_s}\ {\Phi_s(x)}=\frac{1}{2}\|A_sx_s-b\|_2^2+\frac{{\lambda}^{2}}{2}\|x_s\|_2^2%\tag{{$\hat{\mathcal{P}}$}}.        
\end{equation}
{From the numerical examples in Section \ref{sec: Numerical results}, we can see that $x_s$ from \eqref{equivalent vector minimization model} and $x$ from \eqref{vector minimization model} can provide indistinguishable image restoration results. {This is because the original ill-posed problem \eqref{vector minimization model} only requires a numerical solution $x$ with relatively low accuracy. As long as the difference between $A_s$ and $A$ falls below a certain level, which can be easily met with only a small value of $s$ in \eqref{Kronecker structure}, $x_s$ from \eqref{equivalent vector minimization model} and $x$ from \eqref{vector minimization model} can reach the same level of accuracy. This phenomenon is analyzed in Theorem \ref{theorem: total error of sFISTA} in Section \ref{sec: fgFISTA} and verified by the numerical experiments in Section \ref{sec: Numerical results}.}
}

If $b=\mathrm{vec}(B),\ x_s=\mathrm{vec}(X) \ \text{and}\ e=\mathrm{vec}(E)$, where $\mathrm{vec}(X)$ represents a column vector obtained from vectorizing a matrix $X$ ({\rm{i.e.}} columns of $X$ are stacked one after the other), then \eqref{linear system approx} can be rewritten equivalently as
\begin{equation}\label{matrix equation}
B=\sum\limits_{i=1}^{s}H_iXK_i^T+E.        
\end{equation}
It is straightforward to derive the corresponding Tikhonov regularized minimization model as follows
\begin{equation}\label{matrix minimization model}
\min_{X}\ \frac{1}{2}\|\sum\limits_{i=1}^{s}H_iXK_i^T-B\|_F^2+\frac{{\lambda}^{2}}{2}\|X\|_F^2,%\tag{{$\tilde{\mathcal{P}}$}}        
\end{equation}
where $\|\cdot\|_F$ denotes the Frobenius norm. \eqref{matrix minimization model} has several advantages over the original optimization problem \eqref{vector minimization model}. First, \eqref{matrix minimization model} benefits from the Kronecker product structure of $A_s$ and can exploit more computationally intensive matrix-matrix operations. In addition, all the matrices $H_i$ and $K_i$ are structured matrices, which enables fast matrix-vector multiplications at each iteration. Second, the summation of $s$ terms in \eqref{matrix minimization model} can be performed independently and enables \eqref{matrix minimization model} to reach superior parallel efficiency when implemented on modern high performance computing architectures. 
%\JN{Cite some previous work on matrix equations for inverse problems.}
Some work has been done to exploit matrix equation structures for iterative methods to solve inverse problems of the form (\ref{matrix equation}); see, for example,
\cite{Bentbib2018,BoJb07,CaRe96,Zhang2018MatrixEquation}.
{In this paper, we propose the structured FISTA (sFISTA) method. It gains its efficiency by exploiting both the Kronecker product structure of $A$ as well as the structures from $K_i$ and $H_i$. The convergence rate of sFISTA can be of the same order as FISTA under mild conditions.}

The remaining sections are organized as follows. In Section \ref{sec: Kroneker decomposition}, we describe how to approximate the blurring matrix $A$ into the sum of a few of Kronecker products. In Section \ref{sec: fgFISTA}, we first briefly review the FISTA framework and then propose the sFISTA method. We also show that sFISTA for \eqref{matrix minimization model} is equivalent to FISTA for \eqref{equivalent vector minimization model} and derive the convergence and complexity analysis of sFISTA for \eqref{matrix minimization model}. Some numerical examples are provided in Section \ref{sec: Numerical results} and the concluding remarks are drawn in Section \ref{sec:Conclusion}.

\section{Kronecker Decomposition}\label{sec: Kroneker decomposition}

Consider a 2-D spatially-invariant image restoration problem. It was shown in \cite{Ng1999} that three different structures of the blurring matrix $A$ commonly occur. If the zero boundary condition (corresponding to assuming the values of $x$ outside the domain of consideration are zero) is applied, $A$ will be a block-Toeplitz-Toeplitz-block (BTTB) matrix. On the other hand, if the periodic boundary condition (corresponding to the case that the image outside the domain of consideration is a repeat of the image inside in all directions) is used, $A$ becomes a block-circulant-circulant-block (BCCB) matrix. Finally, $A$ would be block-Toeplitz-plus-Hankel with Toeplitz-plus-Hankel-blocks (BTHTHB) if the reflective boundary condition (corresponding to a reflection of the original scene at the boundary) is utilized. In any case, the matrix $A$ can always be approximated as the sum of a few Kronecker products. Since the periodic boundary condition often cause severe ringing artifacts near image borders, only the other two cases are considered in the remaining sections. 

%\JN{The letter P is used for problem and PSF. One is blue, the other italic math font. But still, this might get a bit confusing. Should we change one?}
In practice, the PSF for images with $m\times n$ pixels is often stored as an $m\times n$ array $P$. When $P$ represents the image of a single bright pixel, the process of taking a picture of such an image is equivalent to computing one column of matrix $A$ with column index $t$, where $t$ depends on the location of the point source. Thus, the structure of $A$ is completely determined by that of $P$. More specifically, suppose $P$ has the SVD decomposition $P=U\Sigma V^T$. Let $u_i$ and $v_i$ be the $i$th columns of matrices $U$ and $V$, respectively and $\sigma_1\geq\sigma_2\geq\cdots\geq\sigma_{\min{(m,n)}}$ be the singular values of $P$. It has been shown that $A$ then admits the following Kronecker decomposition \cite{Kamm2000, Nagy2004}
\begin{equation}
A=\sum\limits_{i=1}^{\min{(m,n)}}\ K_i\otimes H_i,
\label{eq:kron1}
\end{equation}
where $K_i$ and $H_i$ are matrices defined based on $u_i$, $v_i$, $\sigma_i$ and boundary conditions. More details on the structure of $K_i$ and $H_i$ will be provided at the end of this section. Because the singular values of $P$ decay quickly in realistic applications, \eqref{eq:kron1} can be further truncated by keeping only the first $s$ terms
\begin{equation}
A \approx A_s=\sum\limits_{i=1}^{s}\ K_i\otimes H_i.
\label{eq:approx}
\end{equation}

{The approximation error introduced in \eqref{eq:approx} has been well studied in \cite{Kamm2000, Nagy2004}. The analysis in \cite{Kamm2000, Nagy2004} shows that the distance between $A$ and $A_s$ is related to the approximation error of a truncated SVD decomposition of a matrix $\bar{P}$, which is summarized in the following theorem for the square PSF case.
\begin{theorem}\label{Kronecker approximation theorem}\cite[Theorem 3.1]{Nagy2004}.
Assume the blurring matrix $A$ is constructed from a PSF $P$ with center $p_{lq}$ located at $(l,q)$, then for both zero boundary condition and reflective boundary condition, we have
\begin{equation}
\left\|A-\sum\limits_{i=1}^{s}\ K_i\otimes H_i\right\|_F=\left\|\bar{P}-\sum\limits_{i=1}^{s}\sigma_iu_i{v_i}^T\right\|_F,
\end{equation}
where $\bar{P}=W_aPW_b$ with $W_a=\mathrm{diag}\left([\ \sqrt{n-l+1}\ \cdots\ \sqrt{n-1}\ \sqrt{n}\ \sqrt{n-1}\ \cdots\ \sqrt{l}\ ]^T\right)$, $W_b=\mathrm{diag}$\\
$\left([\ \sqrt{n-q+1},\ \cdots\ \sqrt{n-1}\ \sqrt{n}\ \sqrt{n-1}\ \cdots\ \sqrt{q}\ ]^T\right)$ for the zero boundary condition case and $\bar{P}=RPR^T$ with $R$ is the Cholesky factor of the symmetric Toeplitz matrix with its first row as $[n,1,0,1,0,1,\cdots]$ for the reflective boundary condition case. Here $\sum\limits_{i=1}^{s}\sigma_iu_i{v_i}^T$ is the summation of the first $s$ terms in the SVD decomposition of $\bar{P}$.
\end{theorem}
}

{Since the singular values of $\bar{P}$ (as well as $P$) decay quickly to zero for most PSFs, Theorem \ref{Kronecker approximation theorem} guarantees that even a small $s$ in \eqref{eq:approx} could lead to very accurate approximation. Numerical experiments in Section \ref{sec: Numerical results} show that taking $s$ as small as $5$ is enough for the image restoration applications under consideration.
}

At the end of this section, let us take a look at the structure of $K_i$ and $H_i$. If the zero boundary condition is used, $K_i$ and $H_i$ have the following Toeplitz structure
\begin{equation}
K_i=\operatorname{toep}(k_i,l)\quad \text{and}\quad H_i=\operatorname{toep}(h_i,q).
\end{equation}
{In the above equations, $k_i=(\sqrt{\sigma_i}u_i)./\operatorname{diag}{(W_a)}$, $h_i=(\sqrt{\sigma_i}v_i)./\operatorname{diag}{(W_b)}$, where $./$ denotes point-wise division.} And $\operatorname{toep}(c,j)$ denotes a banded Toeplitz matrix whose $j$th column is equal to $c$. For example
\begin{equation}
\operatorname{toep}(c,4)=\left[\begin{array}{ccccccccc}
c_4&&c_3&&c_2&&c_1&&0\\
c_5&&c_4&&c_3&&c_2&&c_1\\
0&&c_5&&c_4&&c_3&&c_2\\
0&&0&&c_5&&c_4&&c_3\\
0&&0&&0&&c_5&&c_4\\
\end{array}\right] \qquad \text{with} \qquad c=\left[\begin{array}{c}
c_1\\
c_2\\
c_3\\
c_4\\
c_5\\
\end{array}\right].
\end{equation}
On the other hand, if the reflective boundary condition is applied,  $K_i$ and $H_i$ are equal to the linear combinations of a Toeplitz matrix and a Hankel matrix:
\begin{equation}
\begin{aligned}
K_i=\operatorname{toep}(k_i,l)+\operatorname{hank}(k_i,l)\quad \text{and} \quad
H_i=\operatorname{toep}(h_i,q)+\operatorname{hank}(h_i,q),
\end{aligned}
\end{equation}
{where $k_i=\sqrt{\sigma_i}R^{-1}u_i$, $h_i=\sqrt{\sigma_i}R^{-1}v_i$} and $\operatorname{hank}(c,j)$ denotes a banded Hankel matrix whose first row and last column are defined by $[c_{j+1},\cdots,c_n,0,\cdots,0]$ and $[0,\cdots,0,c_1,\cdots,c_{j-1}]^T$, respectively. For example
\begin{equation}
\operatorname{hank}(c,3)=\left[\begin{array}{ccccccccc}
c_4&&c_5&&0&&0&&0\\
c_5&&0&&0&&0&&0\\
0&&0&&0&&0&&0\\
0&&0&&0&&0&&c_1\\
0&&0&&0&&c_1&&c_2\\
\end{array}\right] \qquad \text{with} \qquad c=\left[\begin{array}{c}
c_1\\
c_2\\
c_3\\
c_4\\
c_5\\
\end{array}\right].
\end{equation}
The sFISTA to be introduced in the next section will benefit from the fast Toeplitz/Hankel matrix-vector product algorithms when multiplying $K_i$ and $H_i$ with vectors at each iteration.

\section{Structured FISTA}\label{sec: fgFISTA}
In this section, we will first review the FISTA framework for solving \eqref{vector minimization model} and then propose sFISTA for solving \eqref{matrix minimization model}. We can prove that the proposed sFISTA for solving \eqref{matrix minimization model} is equivalent to FISTA for solving \eqref{equivalent vector minimization model}. {A detailed error analysis has also been conducted to show that the computational accuracy of sFISTA can reach the same level as that of FISTA under mild conditions. Finally, we compare the computational complexity of sFISTA for solving \eqref{matrix minimization model} and FISTA for solving \eqref{vector minimization model} and show that sFISTA is more efficient in both serial and parallel computing environments.}

\subsection{FISTA: A fast iterative shrinkage-thresholding algorithm}
A fast iterative shrinkage-thresholding algorithm (FISTA) was first proposed in \cite{Beck2009A} to solve the following general nonsmooth convex optimization model
\begin{equation}\label{fista general model}
\min_{x}\ \{F(x)=f(x)+g(x)\},       
\end{equation}
where $f:{\mathbb{R}^N\rightarrow\mathbb{R}}$ is a smooth convex function of the type $\mathrm{C^{1,1}}$ and $g:{\mathbb{R}^N\rightarrow\mathbb{R}}$ is a continuous convex function which is possibly nonsmooth. 
%\begin{equation}
%\begin{aligned}
%p_L(y)&=\arg \min\ \{Q_L(x,y):x\in{\mathbb{R}}^n\},\\
%&=\arg \min_{x}\ \left\{g(x)+\frac{L}{2}\left\Arrowvert x-\left(y-\frac{1}{L}\nabla f(y)\right) \right\Arrowvert^2 \right\}.
%\end{aligned}
%\end{equation}
The basic idea of FISTA is that at each iteration, after getting the current iteration point $x_k$, an additional point $y_{k+1}$ is chosen as the linear combination of the current iteration point $x_k$ and the previous iteration point $x_{k-1}$. The next iteration point $x_{k+1}$ is then set as the unique minimizer $p_{L(f)}(y_{k+1})$ of the quadratic approximation $Q_{L(f)}(x,y_{k+1})$ of $F(x)$ at $y_{k+1}$ with
\begin{equation}
Q_{L(f)}(x,y):=f(y)+\langle x-y,\nabla f(y) \rangle+\frac{L(f)}{2}\|x-y\|_2^2+g(x)
\end{equation} 
and $L(f)$ being the Lipschitz constant of $\nabla f$. For more details about FISTA, one can refer to \cite{Beck2009A}.
%The framework of FISTA for (\ref{fista general model}) can be summarized as follows
%\begin{algorithm}[H]\footnotesize
%\caption{FISTA for (\ref{fista general model})}
%\label{algo: Fista algorithm}
%\leftline{Initialization: Give initial point $y_1=x_0\in{\mathbb{R}^{N}}$ and a Lipschitz constant of $\nabla f: L(f)$. Set $k=1$, $t_1=1$.}
%\begin{description}
%\item[\textbf{Step 1}] Compute $x_{k}$ as follows
%\begin{equation*}
%x_{k}=p_L(y_k).
%\end{equation*}
%\item[\textbf{Step 2}] Compute $t_{k+1}$ as follows
%\begin{equation*}
%t_{k+1}=\frac{1+\sqrt{1+4{t_k}^2}}{2}.
%\end{equation*}
%\item[\textbf{Step 3}] Compute $y_{k+1}$ as follows
%\begin{equation*}
%y_{k+1}=x_{k}+\frac{t_k-1}{t_{k+1}}(x_k-x_{k-1}).
%\end{equation*}
%\item[\textbf{Step 4}] If a termination criterion is met, Stop; else, set $k:=k+1$ and go to Step 1.
%\end{description}
%\end{algorithm}

Obviously, (\ref{vector minimization model}) is a special instance of problem (\ref{fista general model}) if we let $f(x)=\frac{1}{2}\|Ax-b\|_2^2$ and $g(x)=\frac{{\lambda}^{2}}{2}\|x\|_2^2$. In this case, the (smallest) Lipschitz constant of the gradient $\nabla f$ is $L_f=\lambda_{\mathrm{max}}(A^TA)$.  Simple calculations lead to
\begin{equation*}
\begin{aligned}
x_k=p_{L_f}(y_k)&=\arg \min_{x}\ \{Q_{L_f}(x,y_k):x\in{\mathbb{R}}^{mn}\}.\\
&=\arg \min_{x}\ \left\{\langle x,A^T(Ay_k-b) \rangle+\frac{L_f}{2}\|x\|_2^2+L_f\langle x,y_k\rangle+\frac{{\lambda}^{2}}{2}\|x\|_2^2\right\},\\
&=\arg\min_{x}\ \left\{\frac{L_f+{\lambda}^2}{2}\left\Arrowvert x-\frac{1}{L_f+{\lambda}^2}\left(L_fy_k-A^T(Ay_k-b)\right)\right\Arrowvert_2^2\right\},\\
&=\frac{1}{L_f+{\lambda}^2}\left(L_fy_k-A^T(Ay_k-b)\right).
\end{aligned}
\end{equation*}
See Algorithm \ref{algo: Fista for P} for a description of FISTA for \eqref{vector minimization model}. 
\begin{algorithm}[t]\footnotesize
\caption{FISTA for (\ref{vector minimization model})}
\label{algo: Fista for P}
\leftline{Initialization: set initial point $y_1=x_0\in{\mathbb{R}^{N}}$, $ L_f= \lambda_{\max}(A^{T}A)$, $k=1$, $t_1=1$.}
\begin{description}
\item[\textbf{Step 1}] Compute $x_{k}$ as follows
\begin{equation*}
x_{k}=\frac{1}{L_f+{\lambda}^2}\left(L_fy_k-A^T(Ay_k-b)\right).
\end{equation*}
\item[\textbf{Step 2}] Compute $t_{k+1}$ as follows
\begin{equation*}
t_{k+1}=\frac{1+\sqrt{1+4{t_k}^2}}{2}.
\end{equation*}
\item[\textbf{Step 3}] Compute $y_{k+1}$ as follows
\begin{equation*}
y_{k+1}=x_{k}+\frac{t_k-1}{t_{k+1}}(x_k-x_{k-1}).
\end{equation*}
\item[\textbf{Step 4}] If a termination criterion is met, Stop; else, set $k:=k+1$ and go to Step 1.
\end{description}
\end{algorithm}

As can be seen from Algorithm \ref{algo: Fista for P}, the total computational cost of FISTA is dominated by matrix-vector multiplications associated with $A$ and $A^{T}$ at Step $1$. Other steps only involve inexpensive vector and scalar operators. Despite its simplicity, FISTA enjoys a fast global convergence rate, which is summarized in Theorem \ref{thm:converence}.
\begin{theorem}{\rm\cite[Theorem 4.4]{Beck2009A}}
Let $\{x_k\}$, $\{y_k\}$ be generated by FISTA. Then for any $k\geq 1$
\begin{equation*}
F(x_k)-F(x_F^*)\leq \frac{2L(f)\|x_0-x_F^*\|_2^2}{(k+1)^2},
\end{equation*}
where $x_F^*$ is the solution of \eqref{fista general model}.
\label{thm:converence}
\end{theorem}

%\textcolor{red}{Say a few words about the benefits of FISTA over other methods for this problem?}

{It is well known that many first order algorithms have very slow convergence rate. From Theorem \ref{thm:converence}, we can see that FISTA is different from classical first order methods in the sense that it preserves a fast global convergence rate $O(1/k^2)$. That is, in order to obtain a numerical solution $x$ such that $F(x)-F(x_F^*)\leq\epsilon$, the number of iterations required by FISTA is at most $\frac{\sqrt{2L(f)}\|x_0-x_F^*\|_2}{\sqrt{\epsilon}}-1$.   In the next section, we will propose the sFISTA which is more efficient for solving \eqref{matrix minimization model}.}

\subsection{Accelerating FISTA by exploiting structures}
In this section, we will show how to adapt the FISTA framework to solve \eqref{matrix minimization model} by exploiting the two hidden structures.
%In this part, we will give a new alternative for solving the problem \eqref{vector minimization model}. 
We first use a Kronecker product approximation $A_s$ of the coefficient matrix $A$ to introduce problem \eqref{equivalent vector minimization model}, which can be equivalently transformed into a matrix problem \eqref{matrix minimization model}. %and then extend FISTA to solve it. 
%The entire algorithm framework is defined as structured fast iterative shrinkage-thresholding algorithm (sFISTA). Although it seems that both the two terms in the object function of problem \eqref{vector minimization model} are smooth, we still call our proposed algorithm as sFISTA since it can be extended to solve the general nonsmooth optimization model \eqref{fista general model}, which will be illustrated in more details later in this section. Introducing the Kronecker product approximation takes full advantage of the hidden structures of the problem, which enables one more level of parallelism into the computation and makes it possible to utilize BLAS3 operation in the following optimization procedure. Moreover, since all the matrices appearing in the Kronecker product approximation are structured matrices (Toeplitz, Hankel, etc.), exploiting their fast matrix-vector multiplications at each iteration is possible.
Consider the following quadratic approximation of the objective function of \eqref{matrix minimization model} at a given point $Y$:
\begin{equation}
\begin{aligned}
Q_L(X,Y):=&\frac{1}{2}\left\|\sum\limits_{i=1}^{s}H_iYK_i^T-B\right\|_F^2+\left\langle X-Y,\sum\limits_{j=1}^{r}H_j^T\left(\sum\limits_{i=1}^{r}H_iYK_i^T-B\right)K_j \right\rangle_F\\
&+\frac{L}{2}\|X-Y\|_F^2+\frac{{\lambda}^{2}}{2}\|X\|_F^2,
\end{aligned}
\end{equation}
where {$L=\lambda_{\mathrm{max}}(A^TA)$} is the Lipschitz constant of the gradient of the first term in the object function of (\ref{matrix minimization model}). Similar to FISTA, we choose the unique minimizer of the quadratic approximation at point $Y_{k+1}$, which is the linear combination of $X_k$ and $X_{k-1}$, as the new iteration point $X_{k+1}$. Mathematically, we set 
\begin{equation*}
Y_{k+1}=X_{k}+\frac{t_k-1}{t_{k+1}}(X_k-X_{k-1}),
\end{equation*}
where $t_k$ and $t_{k+1}$ are parameters updated in the same way as FISTA to make sFISTA maintain the same convergence rate as FISTA for solving \eqref{matrix minimization model} and compute $X_{k+1}$ as 
\begin{equation*}
\begin{aligned}
X_{k+1}&=p_{L}(Y_{k+1})=\arg \min_{X}\ \{Q_{L}(X,Y_{k+1}):X\in{\mathbb{R}}^{m}\times{\mathbb{R}}^{n}\}.\\
&=\arg \min_{X}\ \left\{\left\langle X,\sum\limits_{j=1}^{r}H_j^T\left(\sum\limits_{i=1}^{r}H_iY_{k+1}K_i^T-B\right)K_j \right\rangle_F+\frac{L}{2}\|X\|_F^2+L\langle X,Y_{k+1}\rangle_F+\frac{{\lambda}^{2}}{2}\|X\|_F^2\right\},\\
&=\arg\min_{X}\ \left\{\frac{L+{\lambda}^2}{2}\left\Arrowvert X-\frac{1}{L+{\lambda}^2}\left(LY_{k+1}-\sum\limits_{j=1}^{r}H_j^T\left(\sum\limits_{i=1}^{r}H_iY_{k+1}K_i^T-B\right)K_j\right)\right\Arrowvert_F^2\right\},\\
&=\frac{1}{L+{\lambda}^2}\left(LY_{k+1}-\sum\limits_{j=1}^{r}H_j^T\left(\sum\limits_{i=1}^{r}H_iY_{k+1}K_i^T-B\right)K_j\right).
\end{aligned}
\end{equation*}
Basic steps of sFISTA for {\eqref{matrix minimization model}} are summarized in Algorithm \ref{algo: GLFista for Ptilde}.
\begin{algorithm}[t]\footnotesize
\caption{sFISTA for {\eqref{matrix minimization model}}}
\label{algo: GLFista for Ptilde}
\leftline{Initialization: Compute a Kronecker product approximation $\sum\limits_{i=1}^{s}\ K_i\otimes H_i$ of the coefficient matrix $A$.}
{Give initial point $Y_1=X_0\in{\mathbb{R}^{m}\times\mathbb{R}^{n}}$ and a Lipschitz constant $L$. Set $k=1$, $t_1=1$, {$b=\mathrm{vec}(B)$}.}
\begin{description}
\item[\textbf{Step 1}] Compute $X_{k}$ as follows
\begin{equation*}
X_{k}=\frac{1}{L+{\lambda}^2}\left(LY_k-\sum\limits_{j=1}^{s}H_j^T\left(\sum\limits_{i=1}^{s}H_iY_kK_i^T-B\right)K_j\right).
\end{equation*}
\item[\textbf{Step 2}] Compute $t_{k+1}$ as follows
\begin{equation*}
t_{k+1}=\frac{1+\sqrt{1+4{t_k}^2}}{2}.
\end{equation*}
\item[\textbf{Step 3}] Compute $Y_{k+1}$ as follows
\begin{equation*}
Y_{k+1}=X_{k}+\frac{t_k-1}{t_{k+1}}(X_k-X_{k-1}).
\end{equation*}
\item[\textbf{Step 4}] If a termination criterion is met, Stop; else, set $k:=k+1$ and go to Step 1.\vspace{3mm}
\item[\textbf{Step 5}] Return $x =\mathrm{vec}(X_k)$
\end{description}
\end{algorithm}

Compared with Algorithm \ref{algo: Fista for P}, there are several major differences between sFISTA and FISTA. First of all, the computational cost of Algorithm \ref{algo: Fista for P} is dominated by matrix-vector multiplications while Algorithm \ref{algo: GLFista for Ptilde} can benefit from more computationally intensive matrix-matrix multiplications. Moreover, since $H_i$ and $K_i$ are all structured matrices (Toeplitz, Hankel, etc.), we can further exploit their fast matrix-vector multiplications at Step $1$ in Algorithm \ref{algo: GLFista for Ptilde}. Second, Algorithm \ref{algo: GLFista for Ptilde} decomposes the computation of $X_k$ as the summation of $s$ terms, which can be computed independently. Therefore, we can easily explore two levels of parallelism at each iteration in Algorithm \ref{algo: GLFista for Ptilde}. The first level corresponds to the structured matrix-vector multiplications with multiple vectors and the second level comes from the summation of $s$ terms. This property enables Algorithm \ref{algo: GLFista for Ptilde} to reach superior parallel performance when implemented on modern high performance architectures. Finally, we can prove that sFISTA for (\ref{matrix minimization model}) is equivalent to FISTA for (\ref{equivalent vector minimization model}), which guarantees the fast convergence.
\begin{theorem}\label{thm: equivalence}
sFISTA for {\eqref{matrix minimization model}} and FISTA for \eqref{equivalent vector minimization model} provide the same output as long as their initial points satisfy $x_0=\mathrm{vec}(X_0)$. Mathematically, suppose $\{X_k\}$, $\{Y_k\}$ are generated by sFISTA for {\eqref{matrix minimization model}} and $\{x_k\}$, $\{y_k\}$ are obtained by FISTA for \eqref{equivalent vector minimization model}, then we have $x_k=\mathrm{vec}(X_k)$ and $y_k=\mathrm{vec}(Y_k)$.
\begin{proof}
To prove the desired results, we first review two important properties of Kronecker products, which will be used in the proof below.
\begin{equation*}
\begin{aligned}
(H\otimes K)\ \mathrm{vec}(Z)&=\mathrm{vec}(KZH^T),\\
(H\otimes K)^T&=H^T\otimes K^T,
\end{aligned}
\end{equation*}
where $H$, $K$ and $Z$ are matrices of appropriate dimensions. Recall the frameworks of two algorithms, to prove they provide the same output, we only have to show that both algorithms are equivalent at Step $1$. Specifically, we just have to prove
\begin{equation}
\mathrm{vec}\left(\sum\limits_{j=1}^{s}H_j^T\left(\sum\limits_{i=1}^{s}H_iY_kK_i^T-B\right)K_j\right)=A_s^T\left(A_sy_k-b\right).
\end{equation}
Utilizing the Kronecker product properties mentioned above, we can get
\begin{equation}
\begin{aligned}
&\mathrm{vec}\left(\sum\limits_{j=1}^{s}H_j^T\left(\sum\limits_{i=1}^{s}H_iY_kK_i^T-B\right)K_j\right),\\
=&\sum\limits_{j=1}^{s}\left\{\mathrm{vec}\left(\sum\limits_{i=1}^{s}H_j^TH_iY_kK_i^TK_j\right)-\mathrm{vec}(H_j^TBK_j)\right\},\\
=&\sum\limits_{j=1}^{s}\left\{\sum\limits_{i=1}^{s}\mathrm{vec}\left(H_j^TH_iY_kK_i^TK_j\right)-\mathrm{vec}(H_j^TBK_j)\right\},\\
=&\sum\limits_{j=1}^{s}\left\{\sum\limits_{i=1}^{s}\left((K_j^TK_i)\otimes(H_j^TH_i)\right)\mathrm{vec}(Y_k)-(K_j^T\otimes H_j^T)\mathrm{vec}(B)\right\},\\
=&\sum\limits_{j=1}^{s}\sum\limits_{i=1}^{s}(K_j^T\otimes H_j^T)(K_i\otimes H_i)\mathrm{vec}(Y_k)-\sum\limits_{j=1}^{s}(K_j^T\otimes H_j^T)\mathrm{vec}(B),\\
=&\left(\sum\limits_{j=1}^{s}K_j^T\otimes H_j^T\right)\left(\sum\limits_{i=1}^{s}K_i\otimes H_i\right)\mathrm{vec}(Y_k)-(\sum\limits_{j=1}^{s}K_j^T\otimes H_j^T)\mathrm{vec}(B),\\
=&A_s^T\left(A_sy_k-b\right),
\end{aligned}
\end{equation}
from which we can derive that $x_k=\mathrm{vec}(X_k)$ and $y_k=\mathrm{vec}(Y_k)$. %$\qed$
\end{proof}
\end{theorem}

It is worth pointing out that sFISTA for {\eqref{matrix minimization model}} is only equivalent to FISTA for \eqref{equivalent vector minimization model} due to the Kronecker product approximation error. {The total computational error of sFISTA for solving \eqref{matrix minimization model} comes from two places: the Kronecker product approximation to $A$ and the iterative procedure of sFISTA. The following theorem analyzes the effect of these two kinds of errors on the accuracy of the final computed result.}

\begin{theorem}\label{theorem: total error of sFISTA}
Assume $x^*$ and $x_s^*$ are the exact solutions of \eqref{vector minimization model} and \eqref{equivalent vector minimization model} respectively, $\{{X}_k,{Y}_k\}$ is the sequence obtained by sFISTA for \eqref{matrix minimization model}. Denote $\tilde{x}_k=\mathrm{vec}(X_k)$. If the singular values of $A^TA_s+\lambda^2I$ have a lower bound and the Kronecker product approximation $A_s$ satisfies $\|A-A_s\|_{F}=\epsilon_s$, then we have for any $k\geq1$
\begin{equation}\label{eqn: theorem4 desired result}
|\Phi(\tilde{x}_k)-\Phi(x^*)|\leq \frac{2L\|\tilde{x}_0-x_s^*\|_2^2}{(k+1)^2}+c_0\epsilon_s,
\end{equation}
where $c_0$ is a positive constant independent of $k$.
\begin{proof}
Since $x^*$ and $x_s^*$ are the exact solutions of \eqref{vector minimization model} and \eqref{equivalent vector minimization model}, respectively, from their optimality conditions we have
\begin{equation}\label{eqn: optimal conditions}
A^T(Ax^*-b)+\lambda^2x^*=0\quad {\rm{and}}\quad A_s^T(A_sx_s^*-b)+\lambda^2x_s^*=0,
\end{equation}
which implies that
\begin{equation}\label{eqn: equivalent optimal conditions}
A^TAx^*=A^Tb-\lambda^2x^*\quad {\rm{and}}\quad A_s^TA_sx_s^*=A_s^Tb-\lambda^2x_s^*.
\end{equation}
It is easy to see
\begin{equation}\label{eqn: Phi decomposition}
\left|\Phi(\tilde{x}_k)-\Phi(x^*)\right|\leq\underbrace{|\Phi(\tilde{x}_k)-\Phi_s(\tilde{x}_k)|}_{I}+\underbrace{|\Phi_s(\tilde{x}_k)-\Phi_s(x_s^*)|}_{II}+\underbrace{|\Phi_s(x_s^*)-\Phi(x^*)|}_{III},\\
\end{equation}

For the first term we have
%\JN{I think we need absolute value here to apply the Cauchy-Schwarz inequality: $|\Phi(\tilde{x}_k-\Phi_s(\tilde{x}_k)|$? Also, to be consistent, subscripts on all norms?}
\begin{equation}
\begin{aligned}
|\Phi(\tilde{x}_k)-\Phi_s(\tilde{x}_k)|&=\left|\frac{1}{2}\|A\tilde{x}_k-b\|_2^2-\frac{1}{2}\|A_s\tilde{x}_k-b\|_2^2\right|,\\
&=\left|\frac{1}{2}(\tilde{x}_k^TA^T-b^T+\tilde{x}_k^TA_s^T-b^T)(A\tilde{x}_k-A_s\tilde{x}_k)\right|,\\
&\leq\left\|\frac{1}{2}\tilde{x}_k^TA^T+\frac{1}{2}\tilde{x}_k^TA_s^T-b^T\right\|_{2}\cdot\|\tilde{x}_k\|_{2}\cdot\|A-A_s\|_2,\\
&\leq c_1\epsilon_s,
\end{aligned}
\end{equation}
where we use the fact that $\|A\|_2$, $\|A_s\|_2$, $\|\tilde{x}_k\|_2$ are bounded.

From Theorem \ref{thm: equivalence} we know that sFISTA for {\eqref{matrix minimization model}} is equivalent to FISTA for \eqref{equivalent vector minimization model}, which implies that the second term satisfies
%\JN{Similar to the previous remark, should we use absolute values around terms?}
\begin{equation}
|\Phi_s(\tilde{x}_k)-\Phi_s(x_s^*)|\leq \frac{2L\|\tilde{x}_0-x_s^*\|_2^2}{(k+1)^2}.
\end{equation}
To estimate the last term, we first prove the following fact. From \eqref{eqn: optimal conditions} we get
\begin{equation*}
\begin{aligned}
\lambda^2(x_s^*-x^*)&=A^T(Ax^*-b)-A_s^T(A_sx_s^*-b),\\
&=(A_s^T-A^T)b+A^TAx^*-A_s^TA_sx_s^*,\\
&=(A_s^T-A^T)b+(A^TAx^*-A^TA_sx^*)+(A^TA_sx^*-A^TA_sx_s^*)+(A^TA_sx_s^*-A_s^TA_sx_s^*),\\
&=(A_s^T-A^T)b+A^T(A-A_s)x^*+A^TA_s(x^*-x_s^*)+(A^T-A_s^T)A_sx_s^*,
\end{aligned}
\end{equation*}
which implies that
\begin{equation*}
(A^TA_s+\lambda^2I)(x_s^*-x^*)=(A_s^T-A^T)b+A^T(A-A_s)x^*+(A^T-A_s^T)A_sx_s^*.
\end{equation*}
Then we have
\begin{equation*}
\begin{aligned}
\|x_s^*-x^*\|_2&=\|(A^TA_s+\lambda^2I)^{-1}\cdot\left((A_s^T-A^T)b+A^T(A-A_s)x^*+(A^T-A_s^T)A_sx_s^*\right))\|_2,\\
&\leq\|(A^TA_s+\lambda^2I)^{-1}\|_2\cdot(\|b\|_2+\|A\|_2\|x^*\|_2+\|A_s\|_2\|x_s^*\|_2)\|A-A_s\|_2,\\
&\leq \tilde{c}\epsilon_s,
\end{aligned}
\end{equation*}
where we utilize the boundedness of $\|A\|_2$, $\|A_s\|_2$, $\|b\|_2$, $\|x^*\|_2$, $\|x_s^*\|_2$ and the assumption that the singular values of $A^TA_s+\lambda^2I$ have a lower bound.

Then for the last term we have
\begin{equation*}
\begin{aligned}
&\quad\ |\Phi_s(x_s^*)-\Phi(x^*)|\\
&=\left|\frac{1}{2}\|A_sx_s^*-b\|_2^2+\frac{\lambda^2}{2}\|x_s^*\|_2^2-\frac{1}{2}\|Ax^*-b\|_2^2-\frac{\lambda^2}{2}\|x^*\|_2^2\right|,\\
&=\left|\frac{1}{2}{x_s^*}^TA_s^TA_sx_s^*-b^TA_sx_s^*+\frac{\lambda^2}{2}\|x_s^*\|_2^2-\frac{1}{2}{x^*}^TA^TAx^*+b^TAx^*-\frac{\lambda^2}{2}\|x^*\|_2^2\right|,\\
&=\left|\frac{1}{2}{x_s^*}^T(A_s^Tb-\lambda^2x_s^*)-b^TA_sx_s^*+\frac{\lambda^2}{2}\|x_s^*\|_2^2-\frac{1}{2}{x^*}^T(A^Tb-\lambda^2x^*)+b^TAx^*-\frac{\lambda^2}{2}\|x^*\|_2^2\right|,\\
&=\left|-\frac{1}{2}b^TA_sx_s^*+\frac{1}{2}b^TAx^*\right|,\\
&=\left|\frac{1}{2}b^T(Ax^*-A_sx_s^*)\right|,\\
&=\left|\frac{1}{2}b^T(Ax^*-Ax_s^*+Ax_s^*-A_sx_s^*)\right|,\\
&\leq\frac{1}{2}\|b\|_2(\|A\|_2\|x^*-x_s^*\|_2+\|A-A_s\|_2\|x_s^*\|_2),\\
&\leq c_2\epsilon_s,
\end{aligned}
\end{equation*}
%\JN{Again, absolute value around this difference in third term, and for the final result?}
where we utilize the boundedness of $\|A\|_2$, $\|b\|_2$, $\|x_s^*\|_2$ and the fact that $\|x_s^*-x^*\|_2\leq\tilde{c}\epsilon_s$.

Based on the above analysis above for the three terms in \eqref{eqn: Phi decomposition}, it follows that
\begin{equation}
|\Phi(\tilde{x}_k)-\Phi(x^*)|\leq \frac{2L\|\tilde{x}_0-x_s^*\|_2^2}{(k+1)^2}+(c_1+c_2)\epsilon_s.
\end{equation}
Let $c_0=c_1+c_2$, then the desired result \eqref{eqn: theorem4 desired result} follows. %$\qed$
\end{proof}
\end{theorem}

Theorem \ref{theorem: total error of sFISTA} shows that the error $|\Phi(\tilde{x}_k)-\Phi(x^*)|$ from sFISTA is bounded by two terms: $\frac{2L\|\tilde{x}_0-x_s^*\|_2^2}{(k+1)^2}$ and $c_0\epsilon_s$. The first term decreases as the iteration proceeds while the second term remains constant during the iteration. In order to let the total error $|\Phi(\tilde{x}_k)-\Phi(x^*)|$ fall below a threshold $\epsilon$, we need to make both terms smaller than $\epsilon$. As discussed before, since only a relatively large $\epsilon$ is necessary in these ill-posed inverse problems, a small $s$ would be enough to guarantee $c_0\epsilon_s<\epsilon$. In this sense, the convergence of sFISTA is dominated by the first term $\frac{2L\|\tilde{x}_0-x_s^*\|_2^2}{(k+1)^2}$ and behaves in a similar way as FISTA.

%
%{Although at first glance, compared with FISTA, the convergence rate of sFISTA is not as good as FISTA because of the Kronecker product approximation error $\epsilon_s$, our numerical results together with Theorem \ref{Kronecker approximation theorem} and Theorem \ref{theorem: total error of sFISTA} show that the Kronecker product approximation error will not worse the convergence rate of sFISTA when $s$ only has to be a very small number ($s=5$ is a choice good enough for all cases in our numerical tests). %In fact, the problem we consider in this paper is ill-posed which results in the precision of numerically solving \eqref{linear system} can not be too high, i.e. $\theta$ will be relatively large. While 
%Since the Kronecker product approximation error $\epsilon_s$ can be measured by the difference between a matrix and the summation of its first $s$ terms of its SVD decomposition, whose singular values decrease very quickly to zero, $\epsilon_s$ will decay quickly as $s$ increases, which means $\epsilon_s$ will be relatively small as long as $s$ is very little. When $s$ is big enough to satisfy $\epsilon_s$ small enough, it derives that the convergence rate of sFISTA is of the same order as FISTA.
%}

\begin{figure}[htbp]
\centering
\subfigure[Singular values of $P$]{
\includegraphics[width=0.45\textwidth]{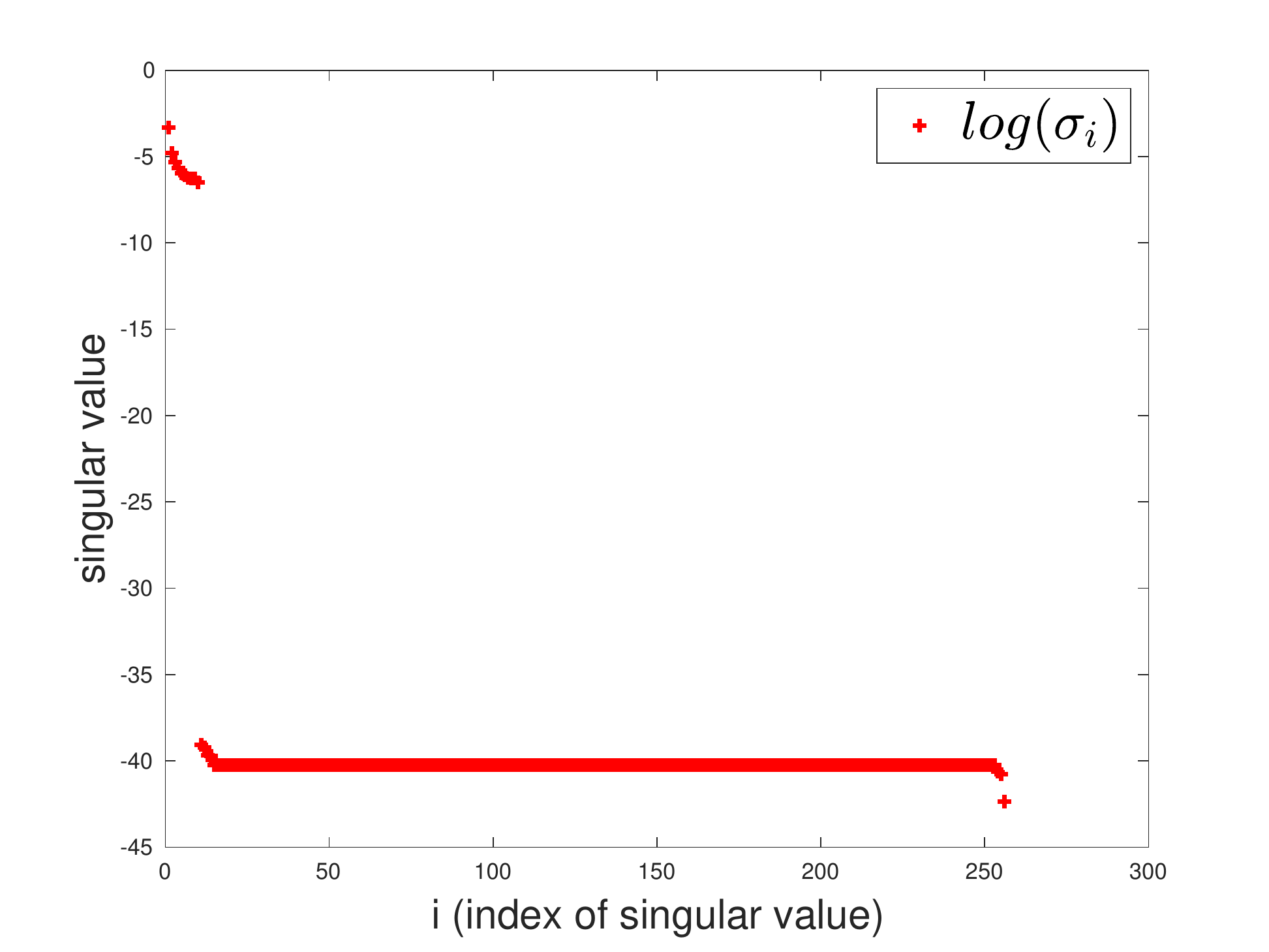}}
\subfigure[Singular values of $\bar{P}$]{
\includegraphics[width=0.45\textwidth]{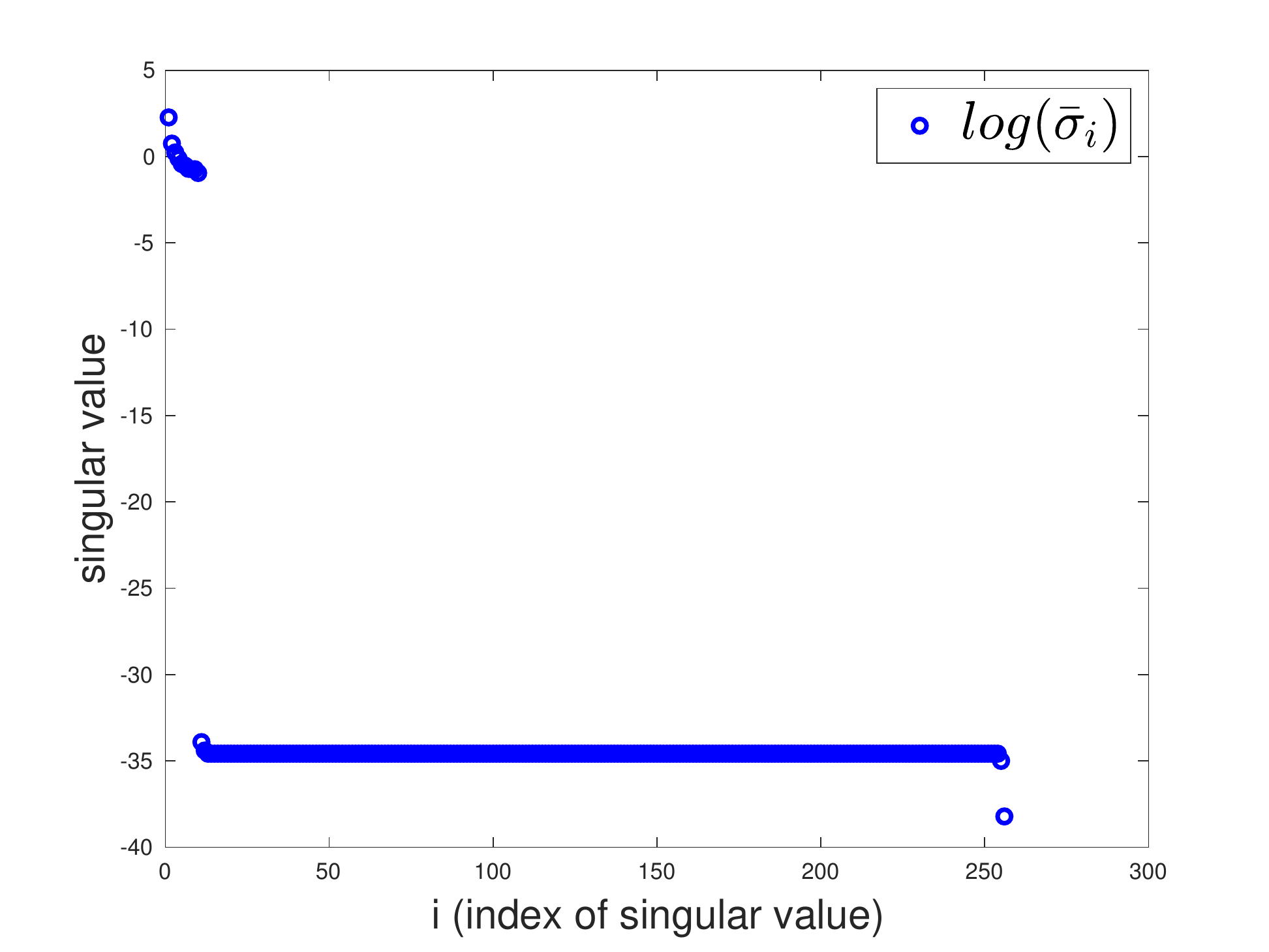}}\\
\caption{Singular values of $P$ and $\bar{P}$ for test image `hst' in Example \ref{example1}.}
\label{figPPbar}
\end{figure}

As an example, we plot the singular values of the matrices $P$ and $\bar{P}$ from the test image `hst' (See Example \ref{example1} in Section \ref{sec: Numerical results} for more details about this image) in Figure \ref{figPPbar}. 
It is easy to see that the singular values of both matrices decay quickly to zero. For example, the ratio of the sixth largest singular value of $P$ to the largest one is only $6.47e\!-\!2$ and the ratio of the tenth largest singular value of $P$ to the largest one reduces to $4.27e\!-\!2$. These patterns can also be observed in other test examples.

\subsection{Complexity Analysis}
In this section, we consider the computational complexity of sFISTA for {\eqref{matrix minimization model}} (Algorithm \ref{algo: GLFista for Ptilde}) and FISTA for (\ref{vector minimization model}) (Algorithm \ref{algo: Fista for P}). If we ignore the structures in $A \in \mathcal{R}^{mn\times mn}$, $K_i\in \mathcal{R}^{n\times n}$, $H_i\in \mathcal{R}^{m\times m}$ and assume that they are all general dense matrices, 
%we can see from the frameworks of Algorithm \ref{algo: %Fista for P} and Algorithm \ref{algo: GLFista for Ptilde} that %the computational complexities of both FISTA and sFISTA %are dominated by step $1$ in the corresponding algorithm. 
then the cost of Step $1$ in Algorithm \ref{algo: Fista for P} and Algorithm \ref{algo: GLFista for Ptilde} would be $O\left(m^2n^2\right)$ and $O\left(s(nm^2+mn^2)\right)$, respectively. When $s$ is much smaller than $m$ and $n$, which is the case for the applications under consideration in this paper, Algorithm \ref{algo: GLFista for Ptilde} is definitely faster than Algorithm \ref{algo: Fista for P}. 
%It is obvious to see that the complexity of sFISTA is better than that of FISTA when $m$ and $n$ are about the same size.

Recall that the blurring matrix $A$ and matrices $K_i$ and $H_i$ from the Kronecker product approximation of $A$ all have specific structures. As the matrix size becomes big enough, these structures will enable us to use fast Fourier transforms (FFTs) to accelerate matrix-vector multiplications encountered in both algorithms. For example, when zero boundary condition is used, $A$ is a block-Toeplitz-Toeplitz-block (BTTB) matrix and $K_i$, $H_i$ are Toeplitz matrices. In this case, the matrix-vector multiplication at step $1$ in Algorithm \ref{algo: Fista for P} can be performed in $O(mn\log(mn))$ with 2D FFTs, while Step $1$ in Algorithm \ref{algo: GLFista for Ptilde} can be done with 1D FFTs in $O\left(smn\log(mn)\right)$. When reflective boundary condition is utilized, $A$ is a block-Toeplitz-plus-Hankel with Toeplitz-plus-Hankel-blocks (BTHTHB) matrix and $K_i$, $H_i$ can be represented as the sum of a Toeplitz matrix and a Hankel matrix. In this case, the computational complexities of Step $1$ in both algorithms are still of the same order as in the zero boundary condition case. Although Algorithm \ref{algo: GLFista for Ptilde} has the same complexity as Algorithm \ref{algo: Fista for P}, it is important to notice that Algorithm \ref{algo: GLFista for Ptilde} is actually much more attractive when implemented on high performance architectures for a number of reasons. First of all, as discussed in the previous section, Algorithm \ref{algo: GLFista for Ptilde} can easily exploit two levels of parallelism, which is crucial for fully taking advantage of the multilelvel parallelism offered by the current architectures. Second, parallel 1D FFTs  are known to scale better than parallel 2D FFTs. Thus, Algorithm \ref{algo: GLFista for Ptilde} is more computationally efficient than Algorithm \ref{algo: Fista for P} for solving large scale problems.
%the sum of $s$ terms in the framework of Algorithm \ref{algo: GLFista for Ptilde} can be computed 
%in a parallel process and also it enables the usage of BLAS3 operation, which can improve the efficiency of Algorithm \ref{algo: GLFista for Ptilde} greatly. The numerical results in Section \ref{sec: Numerical results} verifies that sFISATA is much faster than FISTA.

\section{Numerical Results}
\label{sec: Numerical results}
In this section, we provide some numerical examples to demonstrate the performance of sFISTA for solving {\eqref{matrix minimization model}}.
All the algorithms were implemented with MATLAB and the experiments were performed on a Macbook Air with Intel Core i7 CPU (2.2 GHz). The following notations will be used throughout the section:
\begin{itemize}
\item[$\bullet$] $s$: the number of terms in the Kronecker product approximation;
\item[$\bullet$] $b$: the data vector;
\item[$\bullet$] $\mathrm{noise}$: the vector of perturbations;
\item[$\bullet$] $b_n$: the noisy data $b_n=b+\mathrm{noise}$;
\item[$\bullet$] $\mathrm{NoiseLevel}$: relative level of noise defined as $\|\mathrm{noise}\|_2/\|b\|_2$
\item[$\bullet$] $\mathrm{BlurLevel}$: an indicator used to set the severity of the blur to one of the following: `mild', `medium' and `severe';
\item[$\bullet$] $\eta$: the relative error $\|x-x^*\|/\|x^*\|$;
\item[$\bullet$] $\gamma$: the relative residual $\|r\|/\|b\|$, where $r=Ax-b$;
\item[$\bullet$] $iter$: the iteration number of one algorithm;
\item[$\bullet$] $t(\mathrm{FISTA})$ and $t(\mathrm{sFISTA})$: the CPU time (seconds) of FISTA and sFISTA, respectively;
\item[$\bullet$] $tratio$: an indicator defined as $tratio=\frac{t(\mathrm{sFISTA})}{t(\mathrm{FISTA})}$ to compare the efficiency of FISTA and sFISTA.
\end{itemize}

\begin{exmp}
\label{example1}
In this example, four $256\times 256$ simple test images were extracted based on functions \texttt{PRblurdefocus} and \texttt{PRblurshake} from the regularization toolbox \cite{Gazzola2017I}. The four test images in this example are represented by `hst' (image of the Hubble space telescope), `satellite' (satellite test image), `pattern1' (geometrical image) and `ppower' (random image with patterns of nonzero pixels) respectively, which used reflective (Neumann) boundary conditions \cite{Hansen2006D}. \texttt{PRblurdefocus} and \texttt{PRblurshake} are functions simulating a spatially invariant, out-of-focus blur and spatially invariant motion blur caused by shaking of a camera, respectively. The $\mathrm{BlurLevel}$ was set to be `medium' in these four tests. In addition, function  \texttt{PRnoise} was used to add Gaussian noise with $\mathrm{NoiseLevel}=0.01$ in this example. The regularization parameters were chosen automatically by \texttt{IRhybrid$\_$lsqr} from \cite{Gazzola2017I}, 
which is based on the hybrid bidiagonalization method presented in \cite{HyBR}. 
%\JN{Cite original HyBR paper by Chung, Nagy and O'Leary here.}
The Lipschitz constant was computed as an estimation of the 2-norm of the matrix $A$, which was realized by a few iterations of Lanczos bidiagonalization as implemented in HyBR \cite{Gazzola2017I}. 

We then tested FISTA for (\ref{vector minimization model}) and sFISTA for (\ref{matrix minimization model}) on these four images. To show how the number of terms in the Kronecker product approximation affects the performance of sFISTA, $s$ was set to range from $1$ to $5$ in these four tests. The maximum iteration number for both algorithms was fixed at $50$. To compare the performance of FISTA and sFISTA, we report the CPU time (seconds), the relative error $\eta$ and the relative residual $\gamma$ returned by both algorithms. Their values on these four tests are tabulated in Tables \ref{table1}--\ref{table4}. 

\begin{table}[H]\small
\caption{Numerical results for FISTA and sFISTA for `hst'}
\label{table1}
\begin{center}
\begin{tabular}{@{\extracolsep{\fill}}cccccccccccccc}
\hline
\multirow{2}{*}{} &&\multirow{2}{*}{FISTA} &&\multirow{2}{*}{sFISTA }&&\multirow{2}{*}{sFISTA }&&\multirow{2}{*}{sFISTA }&&\multirow{2}{*}{sFISTA }&&\multirow{2}{*}{sFISTA }\\
&&&&&&&&&&&&\\
&&&&($s=1$)&&($s=2$)&&($s=3$)&&($s=4$)&&($s=5$)\\
\hline
\multirow{2}{*}{time} &&\multirow{2}{*}{$6.8225$}&&\multirow{2}{*}{$0.5255$}&&\multirow{2}{*}{$0.7784$}&&\multirow{2}{*}{$1.0089$}&&\multirow{2}{*}{$1.4913$}&&\multirow{2}{*}{$1.4991$}\\
&&&&&&&&&&&&\\
\multirow{2}{*}{iter} &&\multirow{2}{*}{$50$}&&\multirow{2}{*}{$50$}&&\multirow{2}{*}{$50$}&&\multirow{2}{*}{$50$}&&\multirow{2}{*}{$50$}&&\multirow{2}{*}{$50$}\\
&&&&&&&&&&&&\\
\multirow{2}{*}{$\eta$} &&\multirow{2}{*}{$0.2184$}&&\multirow{2}{*}{$0.2649$}&&\multirow{2}{*}{$0.2329$}&&\multirow{2}{*}{$0.2212$}&&\multirow{2}{*}{$0.2186$}&&\multirow{2}{*}{$0.2182$}\\
&&&&&&&&&&&&\\
\multirow{2}{*}{$\gamma$} &&\multirow{2}{*}{$0.0115$}&&\multirow{2}{*}{$0.0180$}&&\multirow{2}{*}{$0.0123$}&&\multirow{2}{*}{$0.0116$}&&\multirow{2}{*}{$0.0115$}&&\multirow{2}{*}{$0.0115$}\\
&&&&&&&&&&&&\\
\hline
\end{tabular}
\end{center}
\end{table}

\begin{table}[t]\small
\caption{Numerical results for FISTA and sFISTA for `satellite'}
\label{table2}
\begin{center}
\begin{tabular}{@{\extracolsep{\fill}}cccccccccccccc}
\hline
\multirow{2}{*}{} &&\multirow{2}{*}{FISTA} &&\multirow{2}{*}{sFISTA }&&\multirow{2}{*}{sFISTA }&&\multirow{2}{*}{sFISTA }&&\multirow{2}{*}{sFISTA }&&\multirow{2}{*}{sFISTA }\\
&&&&&&&&&&&&\\
&&&&($s=1$)&&($s=2$)&&($s=3$)&&($s=4$)&&($s=5$)\\
\hline
\multirow{2}{*}{time} &&\multirow{2}{*}{$6.8618$}&&\multirow{2}{*}{$0.5593$}&&\multirow{2}{*}{$0.7870$}&&\multirow{2}{*}{$0.9975$}&&\multirow{2}{*}{$1.2309$}&&\multirow{2}{*}{$1.6005$}\\
&&&&&&&&&&&&\\
\multirow{2}{*}{iter} &&\multirow{2}{*}{$50$}&&\multirow{2}{*}{$50$}&&\multirow{2}{*}{$50$}&&\multirow{2}{*}{$50$}&&\multirow{2}{*}{$50$}&&\multirow{2}{*}{$50$}\\
&&&&&&&&&&&&\\
\multirow{2}{*}{$\eta$} &&\multirow{2}{*}{$0.2740$}&&\multirow{2}{*}{$0.3551$}&&\multirow{2}{*}{$0.2979$}&&\multirow{2}{*}{$0.2783$}&&\multirow{2}{*}{$0.2757$}&&\multirow{2}{*}{$0.2744$}\\
&&&&&&&&&&&&\\
\multirow{2}{*}{$\gamma$} &&\multirow{2}{*}{$0.0129$}&&\multirow{2}{*}{$0.0256$}&&\multirow{2}{*}{$0.0152$}&&\multirow{2}{*}{$0.0129$}&&\multirow{2}{*}{$0.0130$}&&\multirow{2}{*}{$0.0129$}\\
&&&&&&&&&&&&\\
\hline
\end{tabular}
\end{center}
\end{table}

\begin{table}[t]\small
\caption{Numerical results for FISTA and sFISTA for `pattern1'}
\label{table3}
\begin{center}
\begin{tabular}{@{\extracolsep{\fill}}cccccccccccccc}
\hline
\multirow{2}{*}{} &&\multirow{2}{*}{FISTA} &&\multirow{2}{*}{sFISTA }&&\multirow{2}{*}{sFISTA }&&\multirow{2}{*}{sFISTA }&&\multirow{2}{*}{sFISTA }&&\multirow{2}{*}{sFISTA }\\
&&&&&&&&&&&&\\
&&&&($s=1$)&&($s=2$)&&($s=3$)&&($s=4$)&&($s=5$)\\
\hline
\multirow{2}{*}{time} &&\multirow{2}{*}{$7.5562$}&&\multirow{2}{*}{$0.6701$}&&\multirow{2}{*}{$0.8669$}&&\multirow{2}{*}{$1.0097$}&&\multirow{2}{*}{$1.2027$}&&\multirow{2}{*}{$1.5385$}\\
&&&&&&&&&&&&\\
\multirow{2}{*}{iter} &&\multirow{2}{*}{$50$}&&\multirow{2}{*}{$50$}&&\multirow{2}{*}{$50$}&&\multirow{2}{*}{$50$}&&\multirow{2}{*}{$50$}&&\multirow{2}{*}{$50$}\\
&&&&&&&&&&&&\\
\multirow{2}{*}{$\eta$} &&\multirow{2}{*}{$0.0607$}&&\multirow{2}{*}{$0.4781$}&&\multirow{2}{*}{$0.2171$}&&\multirow{2}{*}{$0.1417$}&&\multirow{2}{*}{$0.0813$}&&\multirow{2}{*}{$0.0689$}\\
&&&&&&&&&&&&\\
\multirow{2}{*}{$\gamma$} &&\multirow{2}{*}{$0.0087$}&&\multirow{2}{*}{$0.0336$}&&\multirow{2}{*}{$0.0137$}&&\multirow{2}{*}{$0.0135$}&&\multirow{2}{*}{$0.0108$}&&\multirow{2}{*}{$0.0089$}\\
&&&&&&&&&&&&\\
\hline
\end{tabular}
\end{center}
\end{table}

\begin{table}[t]\small
\caption{Numerical results for FISTA and sFISTA for `ppower'}
\label{table4}
\begin{center}
\begin{tabular}{@{\extracolsep{\fill}}cccccccccccccc}
\hline
\multirow{2}{*}{} &&\multirow{2}{*}{FISTA} &&\multirow{2}{*}{sFISTA }&&\multirow{2}{*}{sFISTA }&&\multirow{2}{*}{sFISTA }&&\multirow{2}{*}{sFISTA }&&\multirow{2}{*}{sFISTA }\\
&&&&&&&&&&&&\\
&&&&($s=1$)&&($s=2$)&&($s=3$)&&($s=4$)&&($s=5$)\\
\hline
\multirow{2}{*}{time} &&\multirow{2}{*}{$7.1086$}&&\multirow{2}{*}{$0.5863$}&&\multirow{2}{*}{$0.8070$}&&\multirow{2}{*}{$1.0227$}&&\multirow{2}{*}{$1.2676$}&&\multirow{2}{*}{$1.4446$}\\
&&&&&&&&&&&&\\
\multirow{2}{*}{iter} &&\multirow{2}{*}{$50$}&&\multirow{2}{*}{$50$}&&\multirow{2}{*}{$50$}&&\multirow{2}{*}{$50$}&&\multirow{2}{*}{$50$}&&\multirow{2}{*}{$50$}\\
&&&&&&&&&&&&\\
\multirow{2}{*}{$\eta$} &&\multirow{2}{*}{$0.0968$}&&\multirow{2}{*}{$0.2725$}&&\multirow{2}{*}{$0.1478$}&&\multirow{2}{*}{$0.1160$}&&\multirow{2}{*}{$0.1097$}&&\multirow{2}{*}{$0.0985$}\\
&&&&&&&&&&&&\\
\multirow{2}{*}{$\gamma$} &&\multirow{2}{*}{$0.0096$}&&\multirow{2}{*}{$0.0226$}&&\multirow{2}{*}{$0.0207$}&&\multirow{2}{*}{$0.0139$}&&\multirow{2}{*}{$0.0132$}&&\multirow{2}{*}{$0.0096$}\\
&&&&&&&&&&&&\\
\hline
\end{tabular}
\end{center}
\end{table}

As can be seen from Tables \ref{table1}--\ref{table4},  sFISTA is much faster than FISTA in all test problems. As $s$ increases from $1$ to $5$, the computational time of sFISTA increases monotonically while both the relative error $\eta$ and the relative residual $\gamma$ keep decreasing. When $s$ reaches $5$, the errors of sFISTA are close enough to those of FISTA and sFISTA is still about $5$ times faster than FISTA. We would like to emphasize that sFISTA was only implemented as a serial code and we expect to see a larger speedup with a parallel implementation in the future.

We also plot the four images obtained by sFISTA in Figures \ref{fig1}-\ref{fig4}. As a comparison, the true, blurred and noisy image and the image obtained by FISTA are also provided.  It is easy to see that the images obtained by FISTA and sFISTA ($s=5$) seem very similar to each other.

 \begin{figure}[htbp]
 \centering
 \subfigure[True Image]{
 \includegraphics[width=0.23\textwidth]{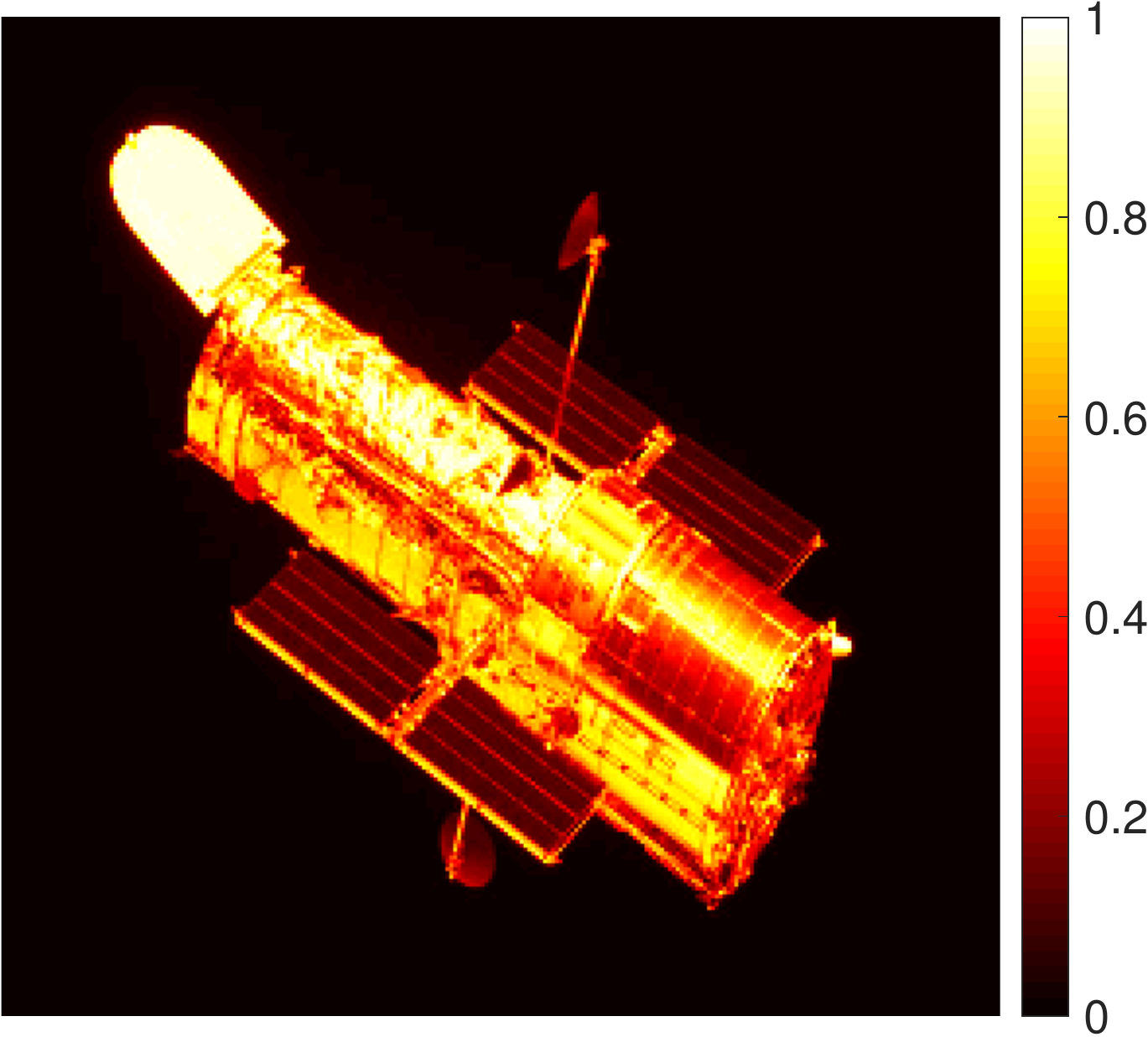}}
 \subfigure[Blurred and noisy Image]{
 \includegraphics[width=0.23\textwidth]{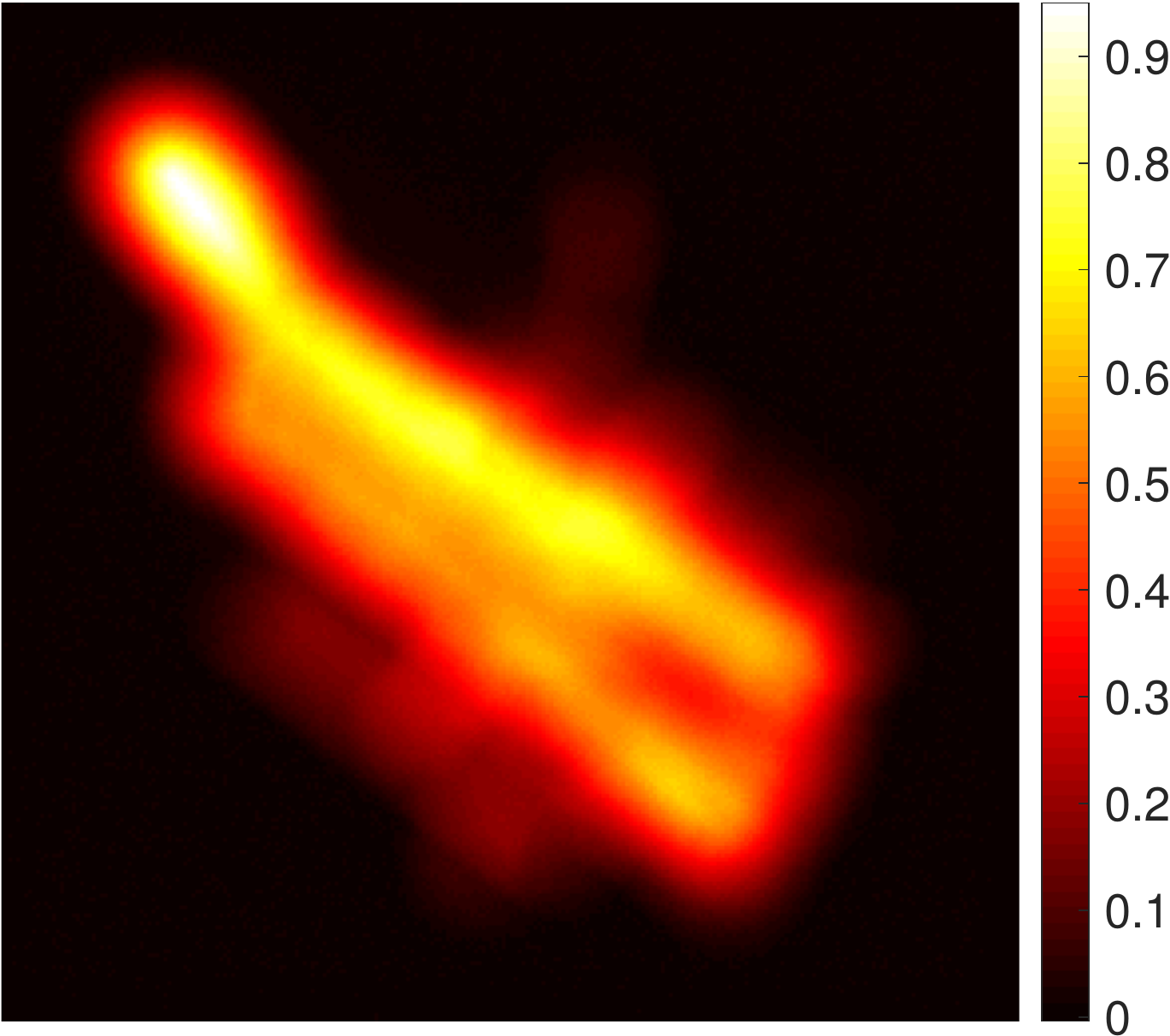}}\\
 \subfigure[Image obtained by FISTA]{
 \includegraphics[width=0.23\textwidth]{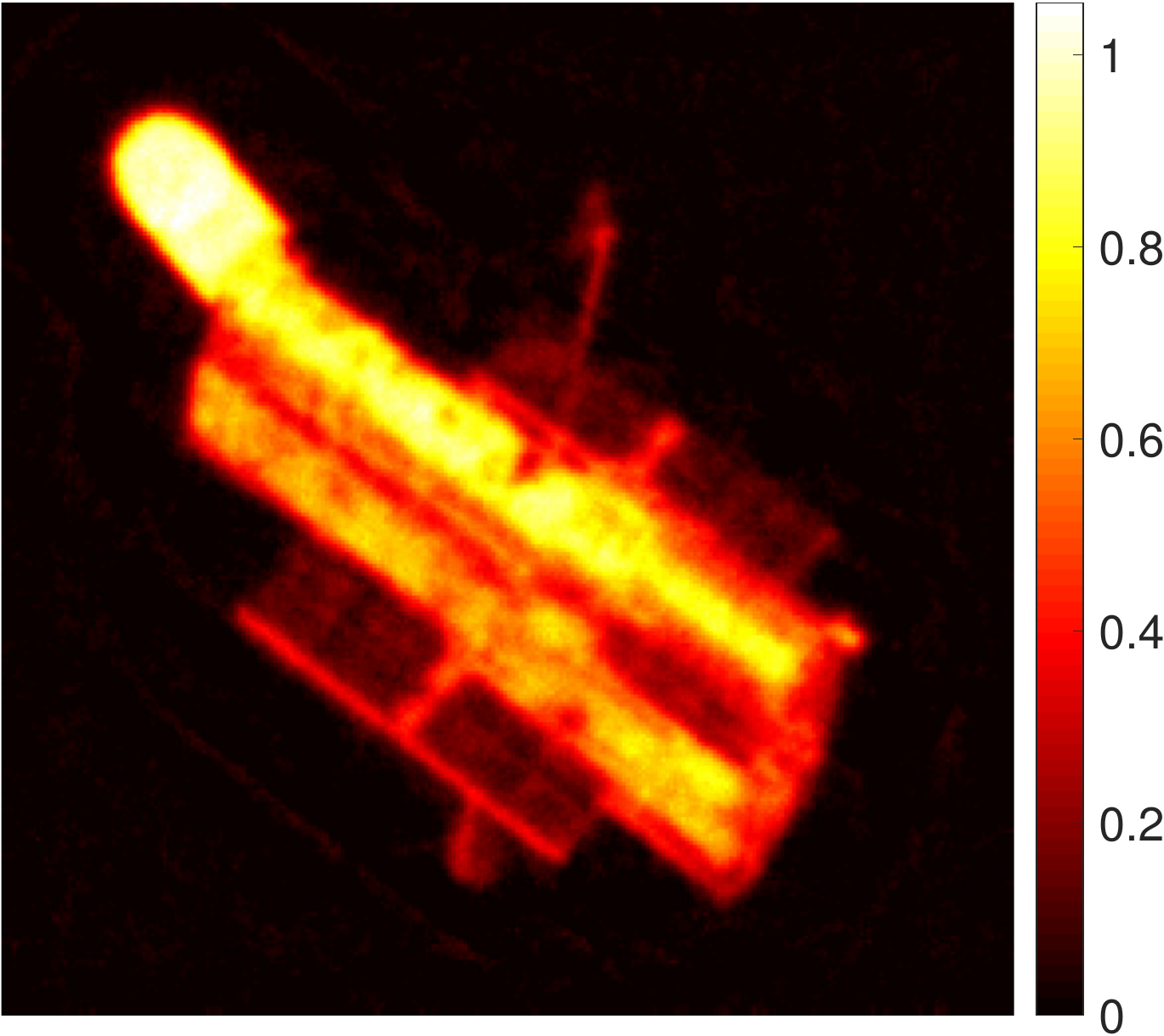}}
 \subfigure[Image obtained by sFISTA  ($s=5$)]{
 \includegraphics[width=0.23\textwidth]{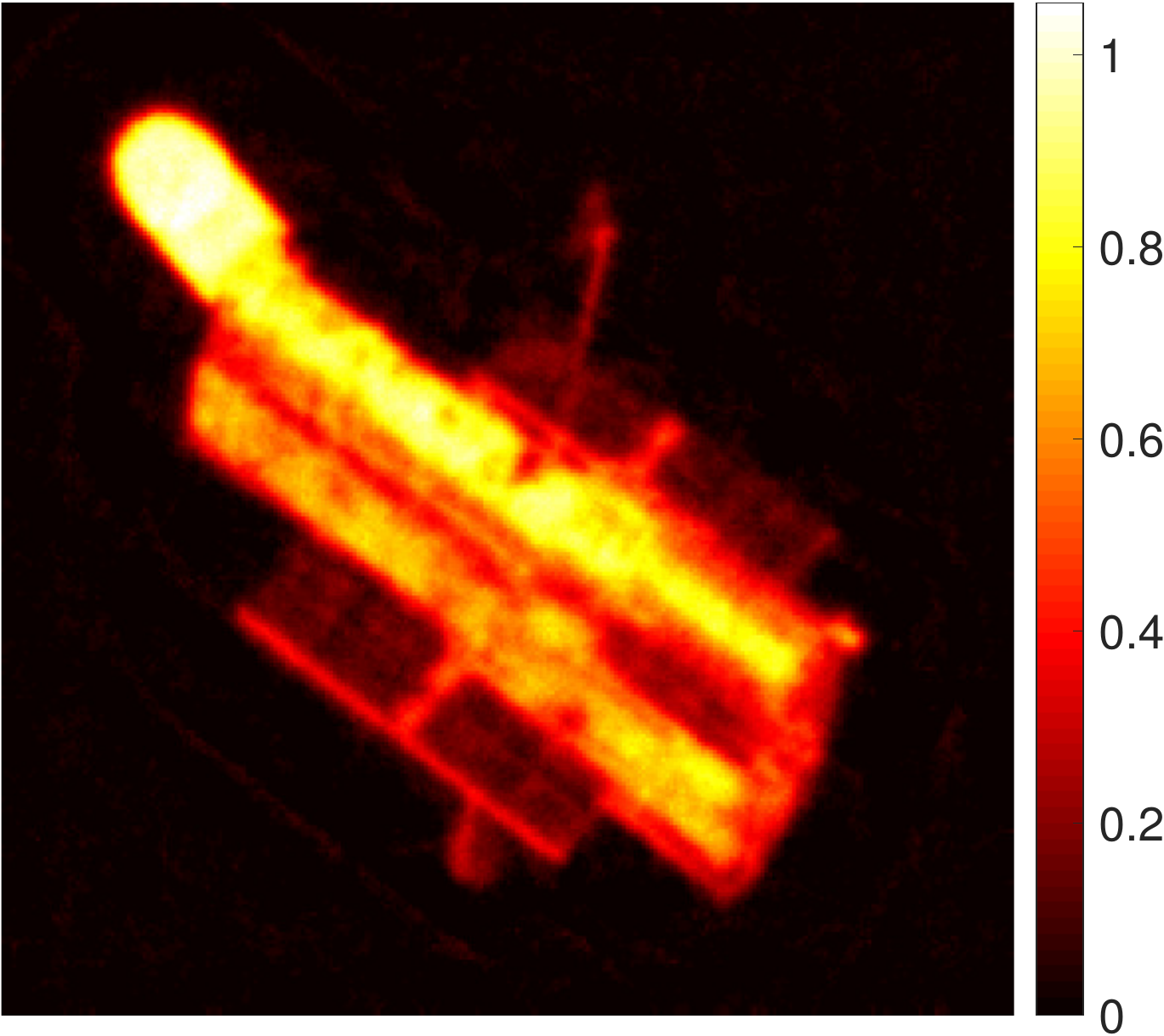}}
 \caption{Figures for `hst' extracted from PRblurdefocus.\label{fig1}}
 \end{figure}

 \begin{figure}[htbp]
 \centering
 \subfigure[True Image]{
 \includegraphics[width=0.23\textwidth]{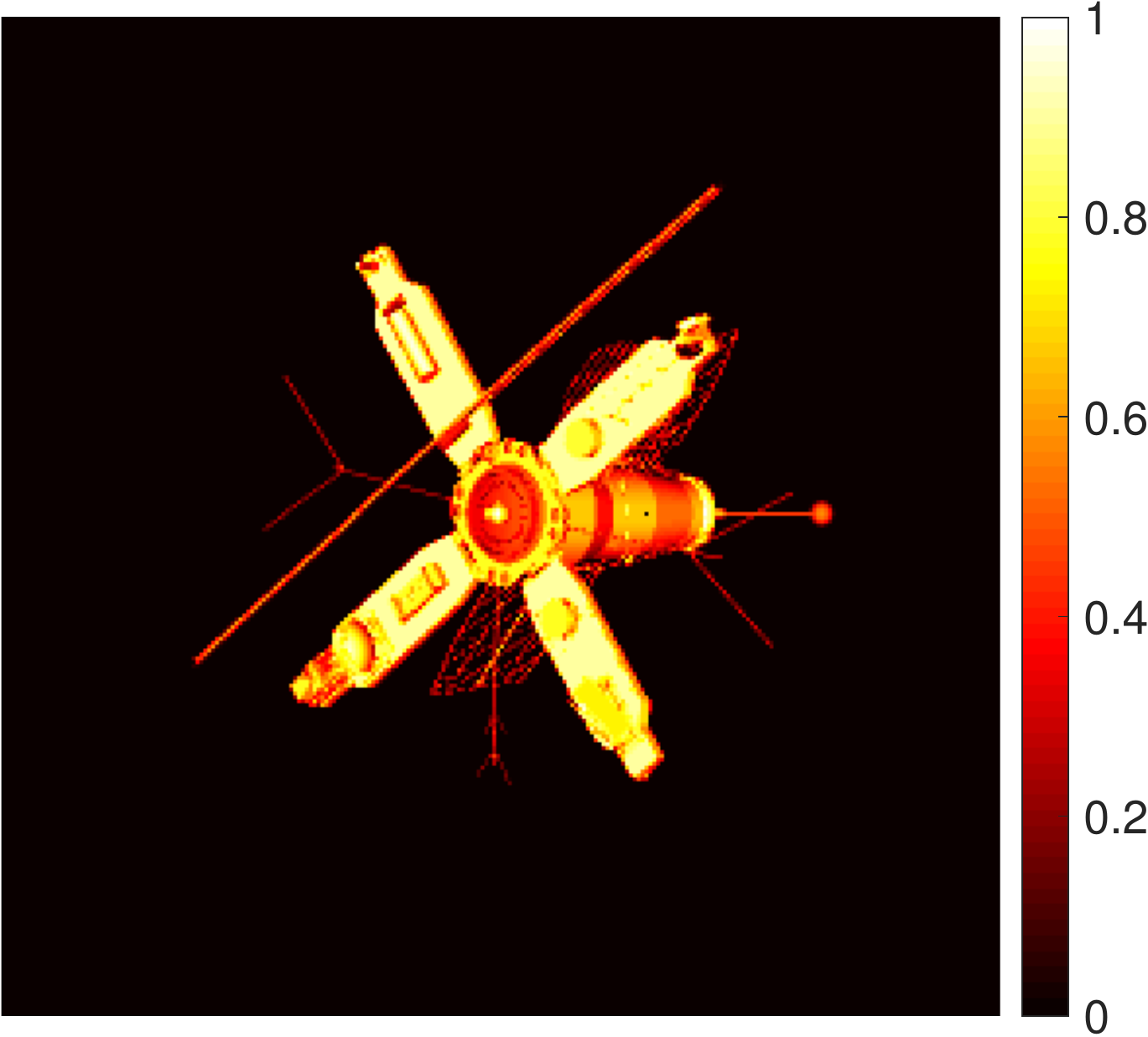}}
 \subfigure[Blurred and noisy Image]{
 \includegraphics[width=0.23\textwidth]{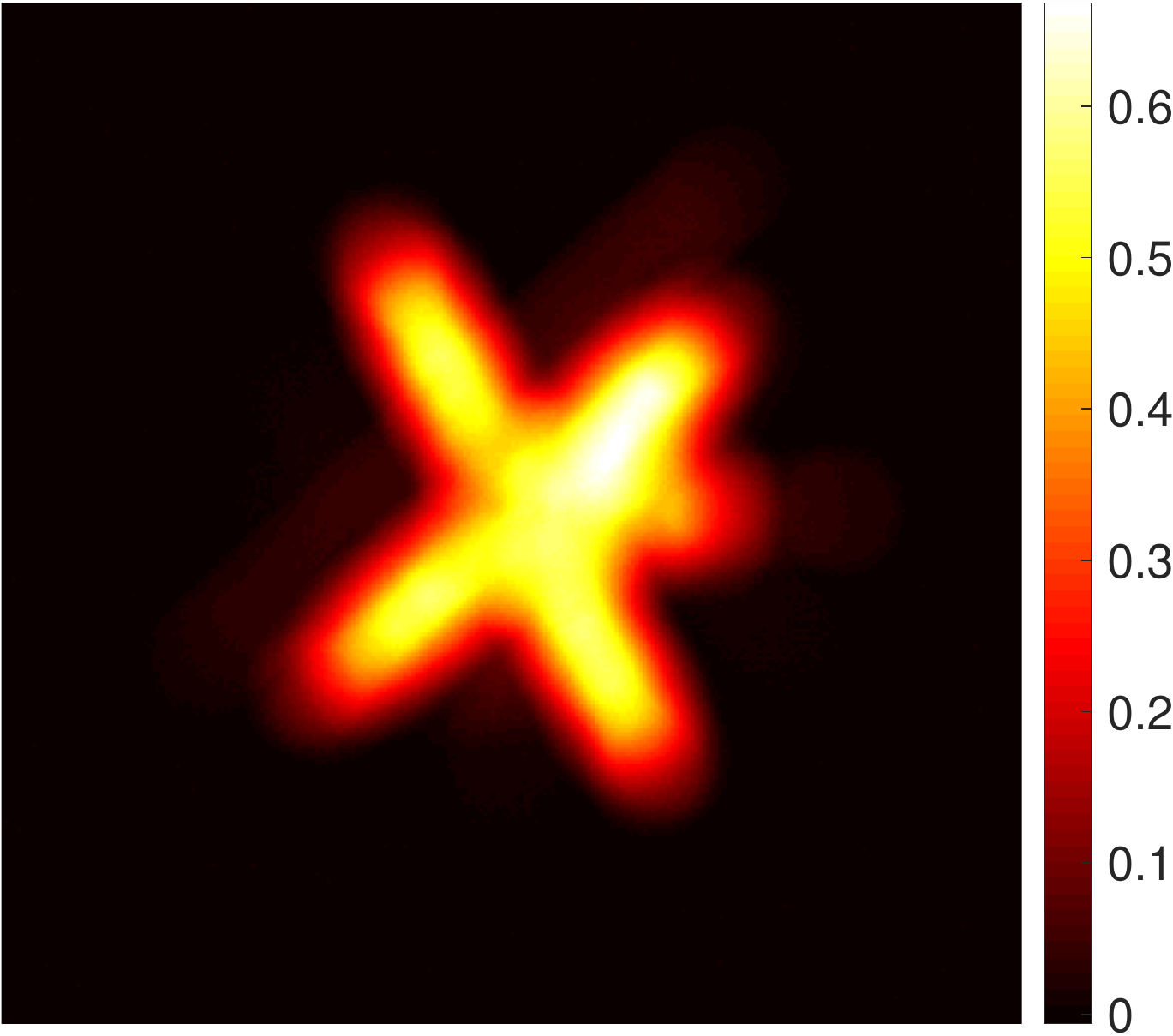}}\\
 \subfigure[Image obtained by FISTA]{
 \includegraphics[width=0.23\textwidth]{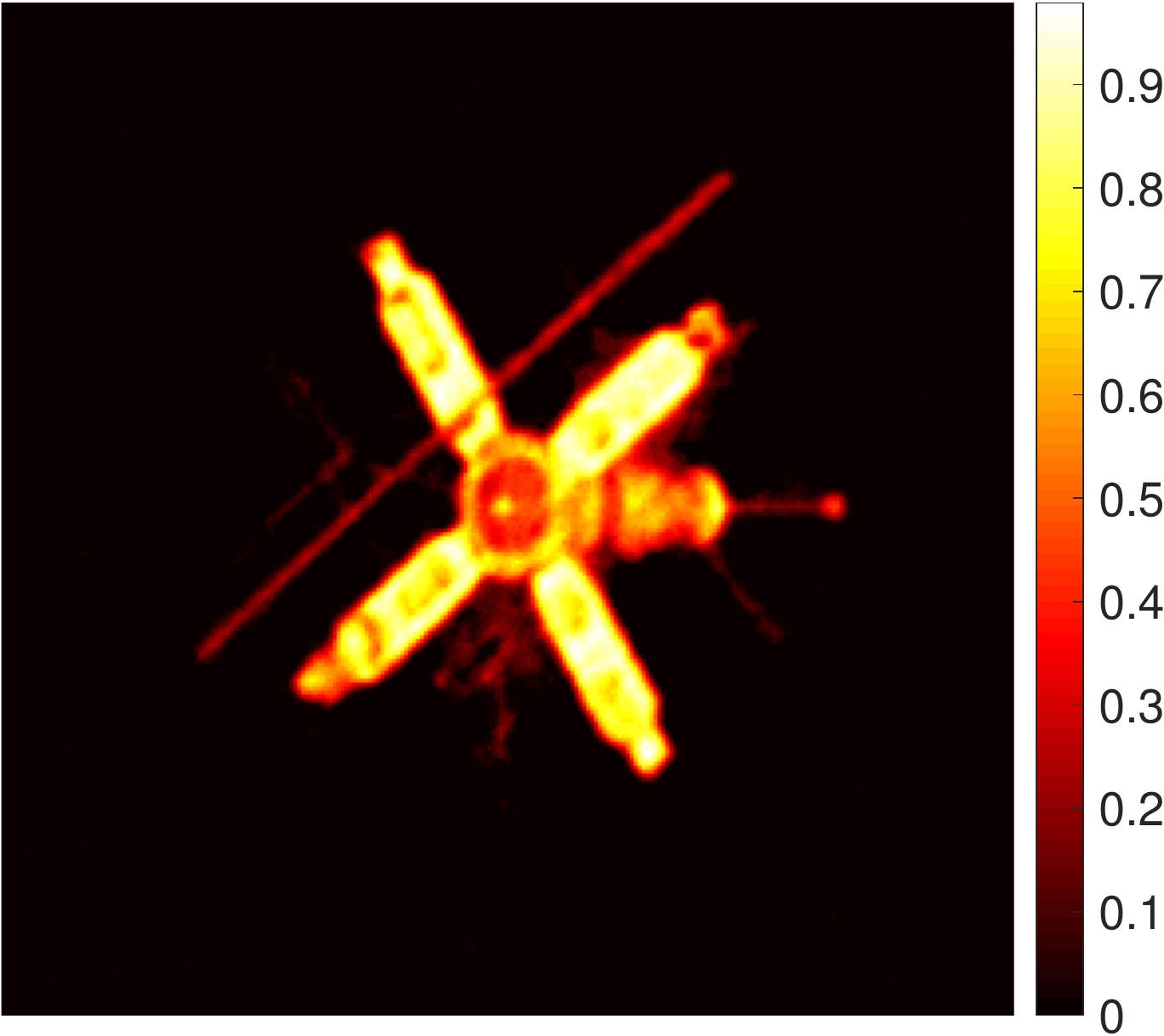}}
 \subfigure[Image obtained by sFISTA ($s=5$)]{
 \includegraphics[width=0.23\textwidth]{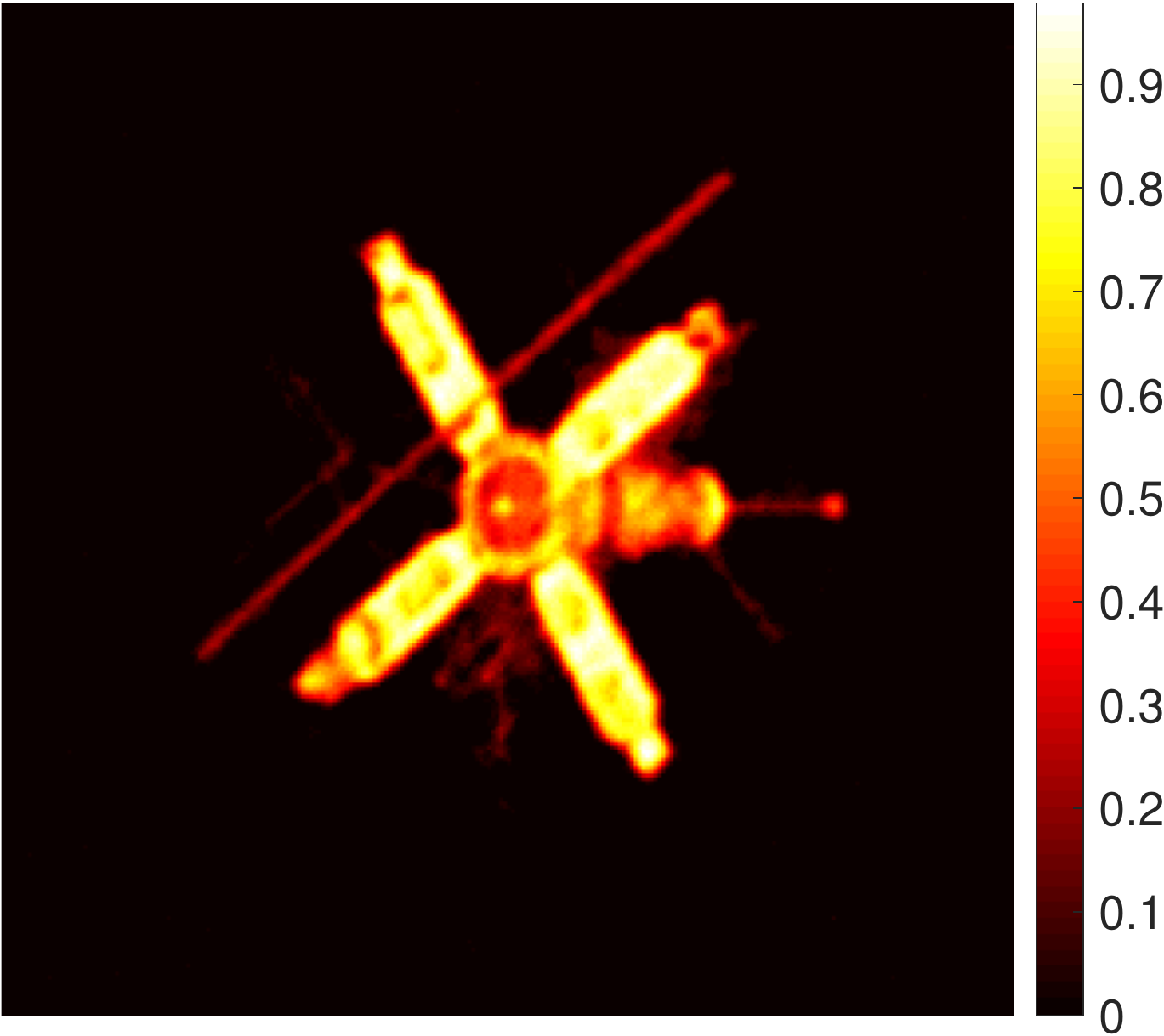}}
 \caption{Figures for `satellite' extracted from PRblurdefocus}
 \label{fig2}
 \end{figure}

 \begin{figure}[htbp]
 \centering
 \subfigure[True Image]{
 \includegraphics[width=0.23\textwidth]{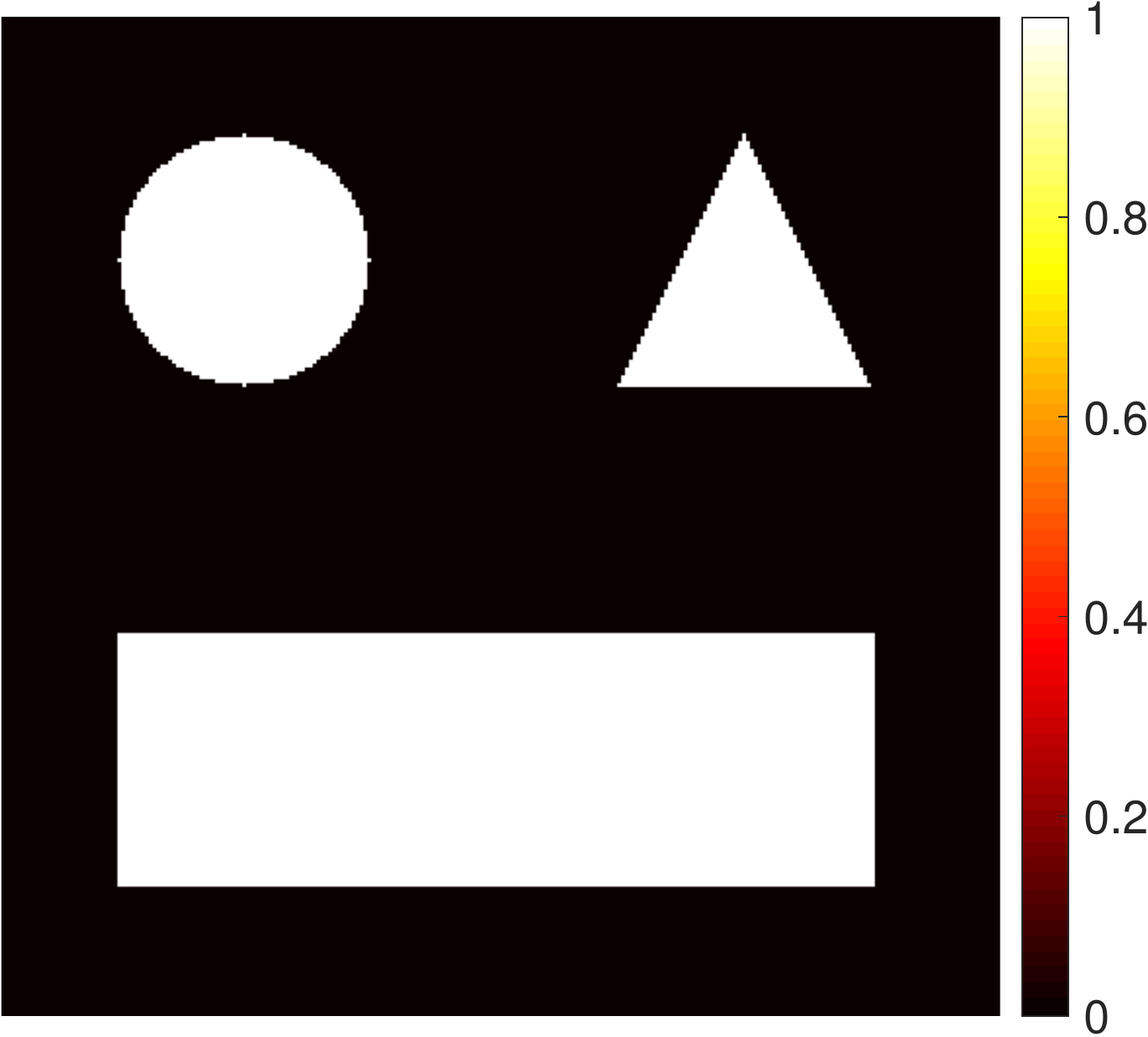}}
 \subfigure[Blurred and noisy Image]{
 \includegraphics[width=0.23\textwidth]{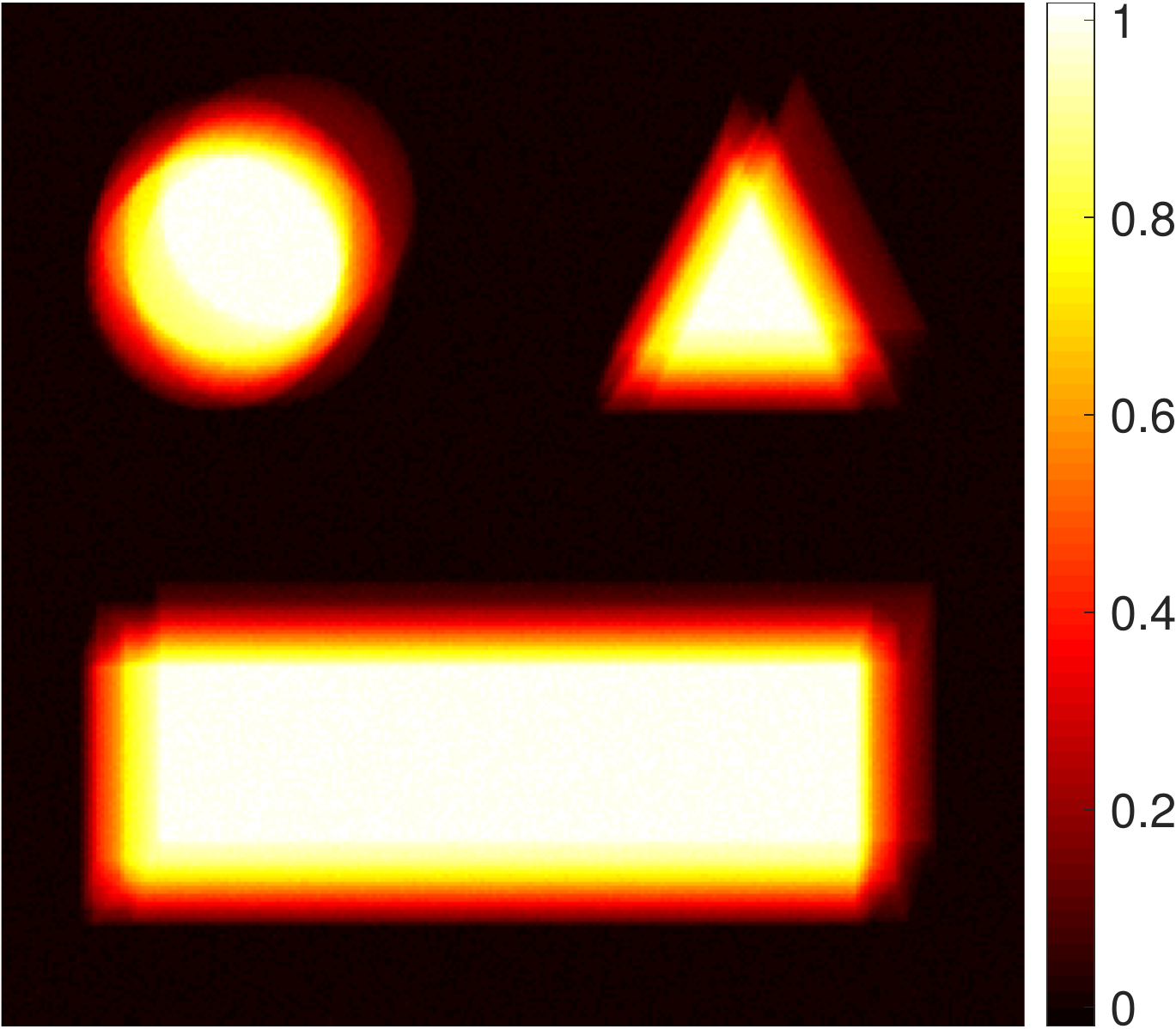}}\\
 \subfigure[Image obtained by FISTA]{
 \includegraphics[width=0.23\textwidth]{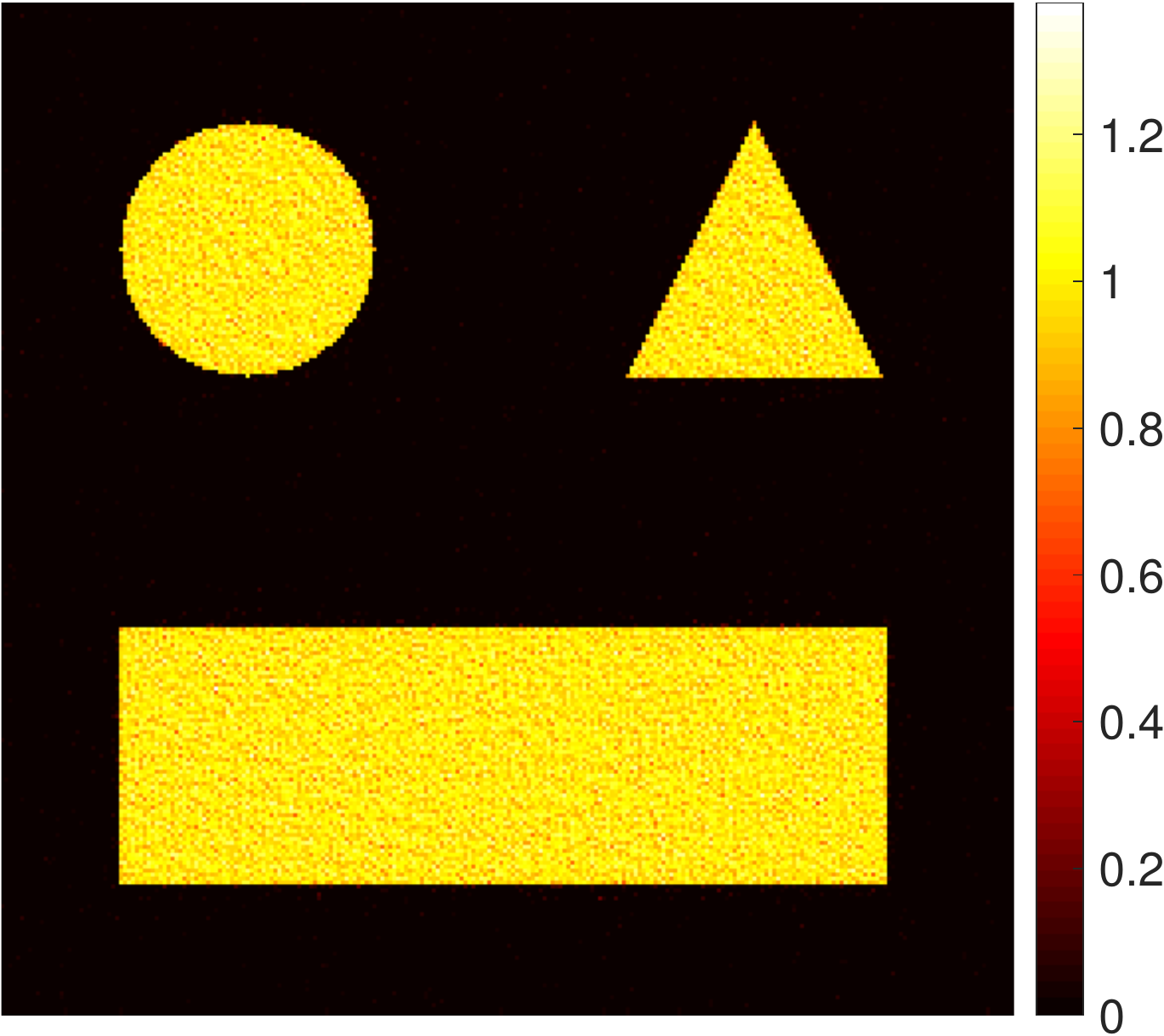}}
 \subfigure[Image obtained by sFISTA ($s=5$)]{
 \includegraphics[width=0.23\textwidth]{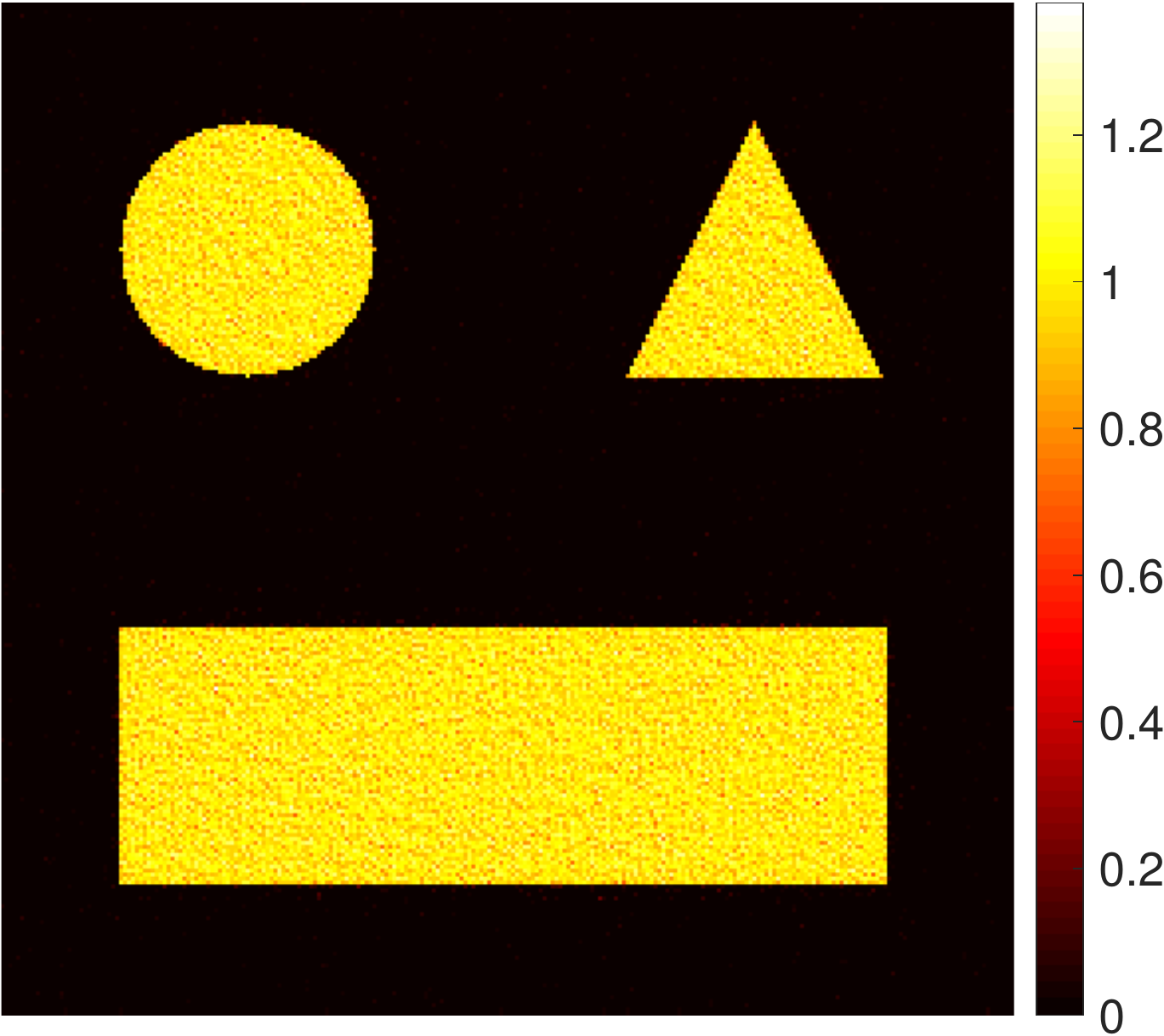}}
 \caption{Figures for `pattern1' extracted from PRblurshake}
 \label{fig3}
 \end{figure}
 
 \begin{figure}[htbp]
 \centering
 \subfigure[True Image]{
 \includegraphics[width=0.23\textwidth]{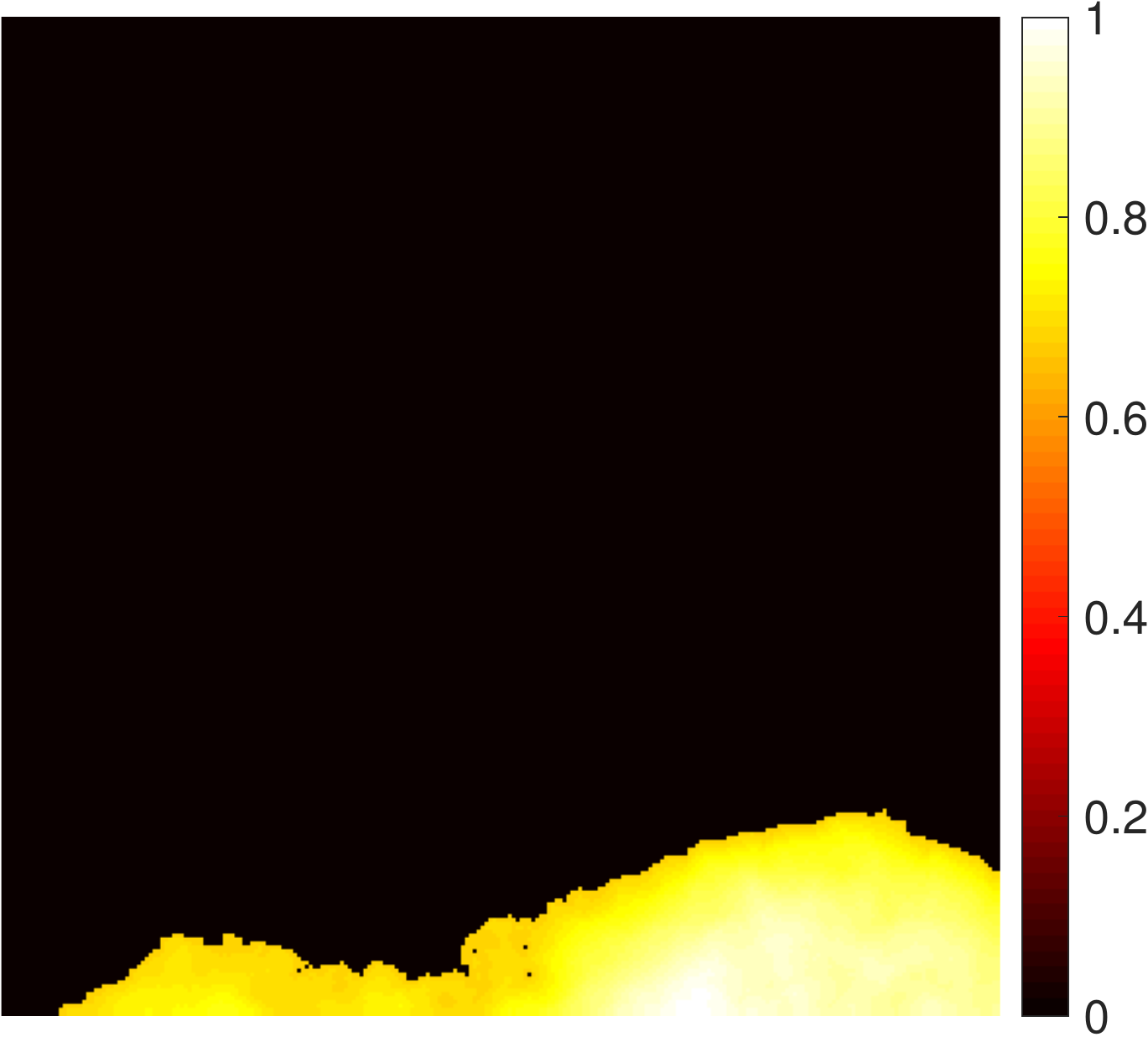}}
 \subfigure[Blurred and noisy Image]{
 \includegraphics[width=0.23\textwidth]{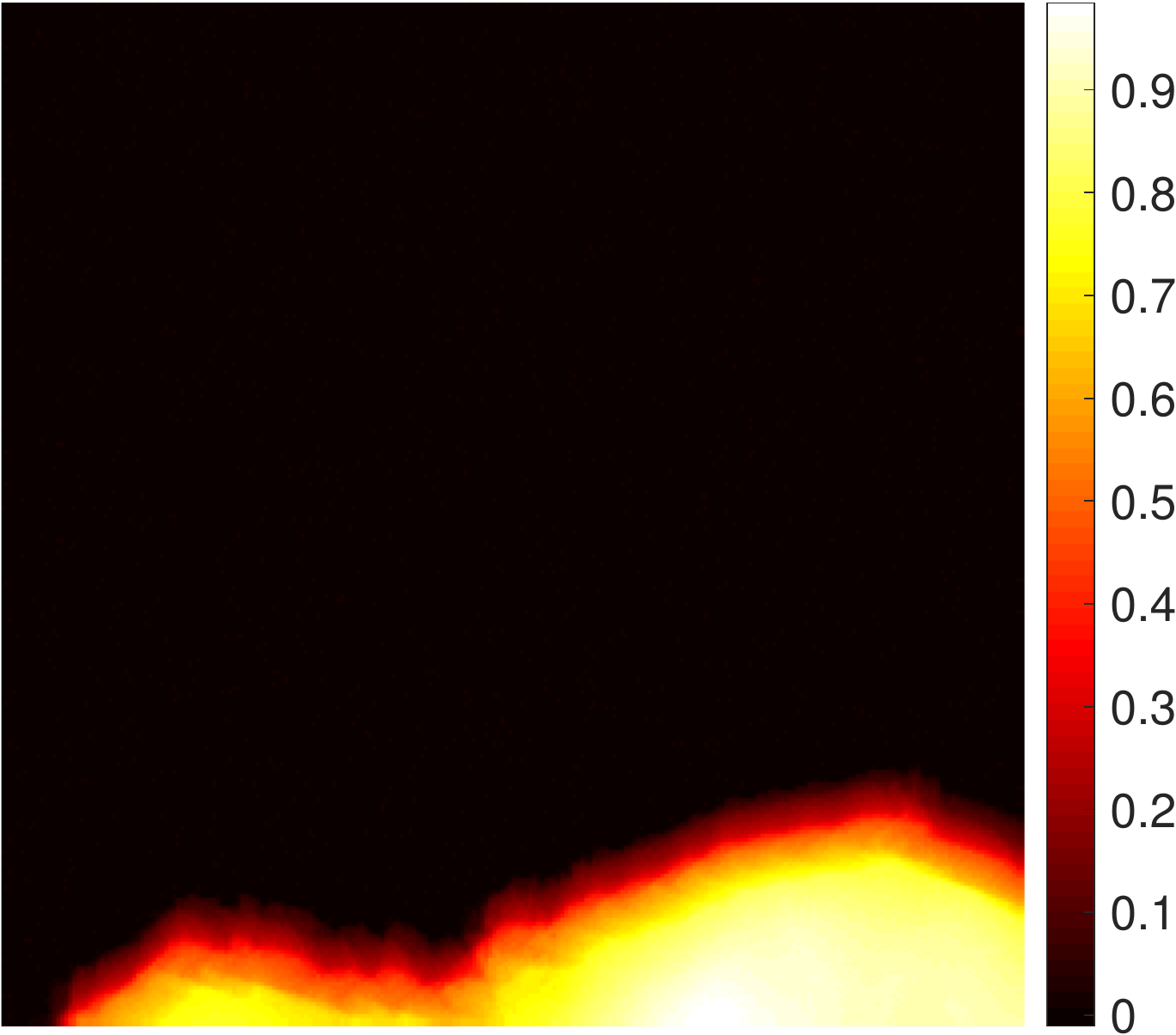}}\\
 \subfigure[Image obtained by FISTA]{
 \includegraphics[width=0.23\textwidth]{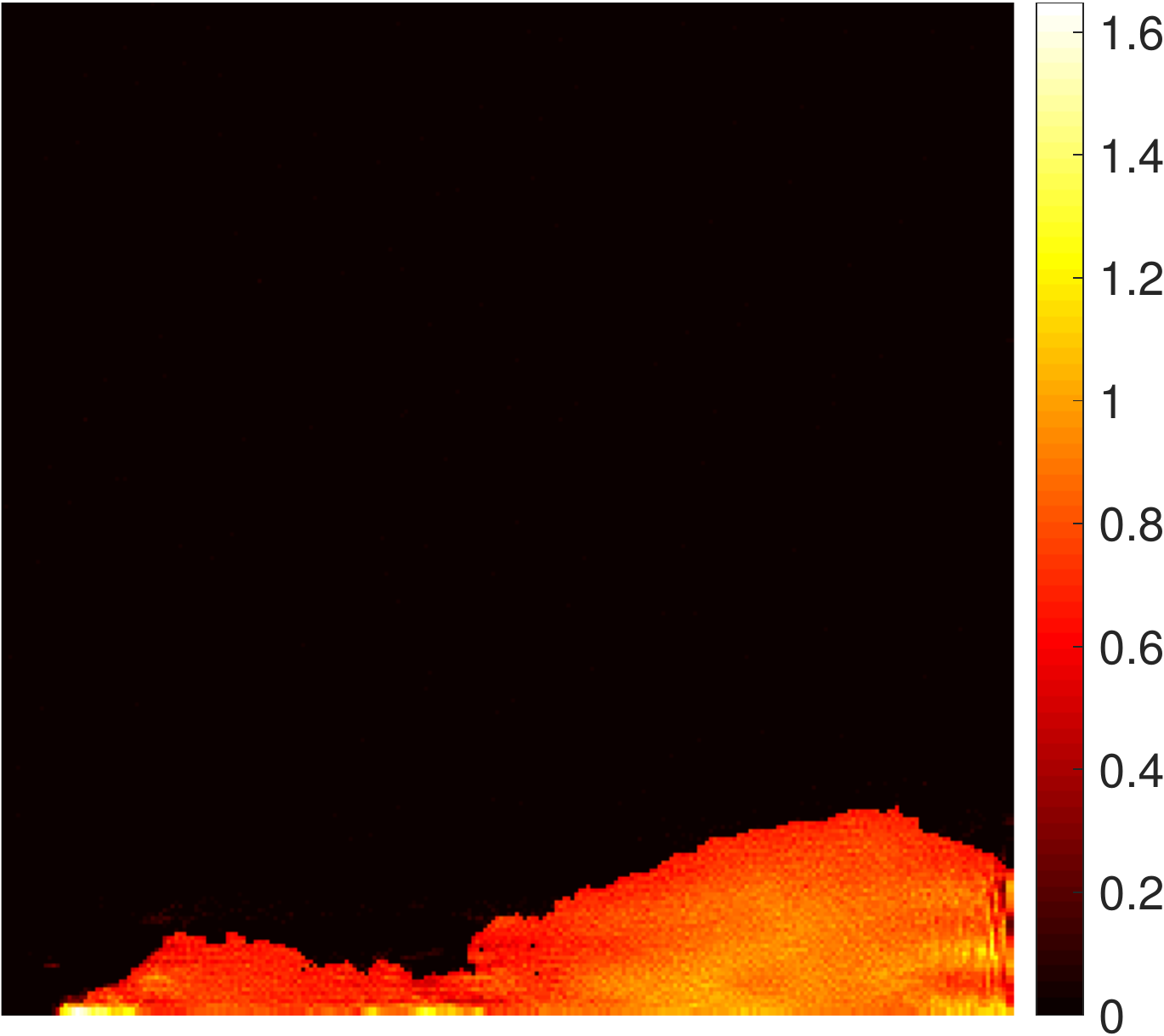}}
 \subfigure[Image obtained by sFISTA ($s=5$)]{
 \includegraphics[width=0.23\textwidth]{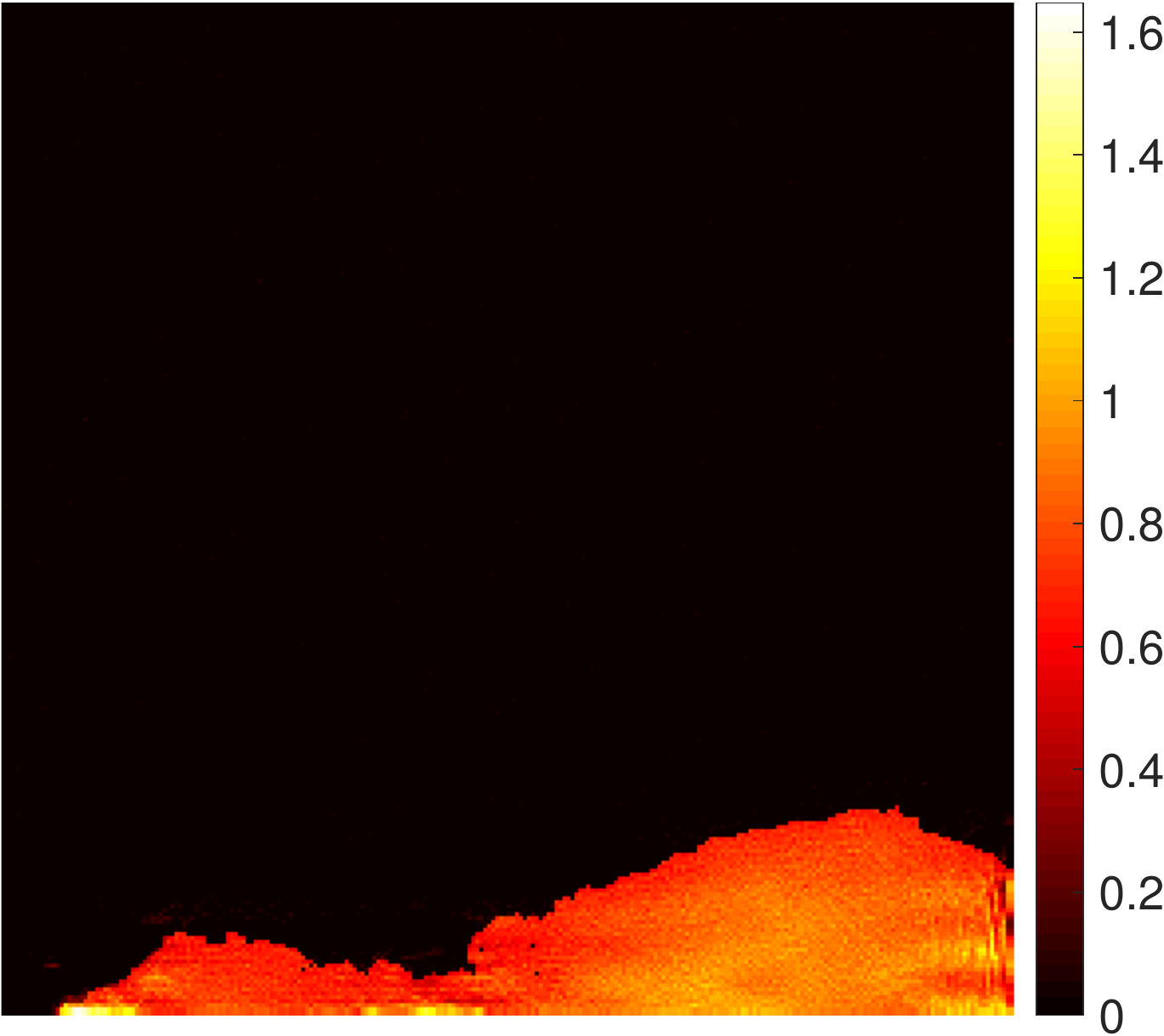}}
 \caption{Figures for `ppower' extracted from PRblurshake}
 \label{fig4}
 \end{figure}

\end{exmp}

\begin{exmp}
In this set of test problems, we compare the performance of sFISTA  and FISTA on eight images with different blur levels and noise levels. These problems are all extracted from functions \texttt{PRblurdefocus} and \texttt{PRblurshake}. The eight test images are represented by `pattern1' (geometrical image), `pattern2' (geometrical image), `ppower' (random image with patterns of nonzero pixels), `smooth' (very smooth image), `dot2' (two small Gaussian shaped dots), `dotk' ($N/2$ small Gaussian shaped dots), `satellite' (satellite test image) and `hst' (image of the Hubble space telescope), respectively. Each test image in this example undergoes blurring and noise-adding procedure with three different blurring levels: `mild', `medium', `severe' and three different kinds of noise `gauss' (Gaussian white noise), `laplace' (Laplacian noise), `multiplicative' (specific type of multiplicative noise, which is often encountered in radar and ultrasound imaging). Since we have seen from Example \ref{example1} that $s=5$ is a good choice for sFISTA, $s$ is fixed to be $5$ in this example. 

We then tested sFISTA and FISTA on these eight images to compare their performance.  The computational results are tabulated in Table \ref{table5} and Table \ref{table6}, respectively. In both tables, there are $9$ cases corresponding to three different $\mathrm{Blurlevel}$ and three different types of noise. Since the indicator $tratio$ measures the ratio of the computational time of sFISTA to that of FISTA for solving the same problem to the same accuracy. Therefore, the smaller $tratio$ is, the more efficient sFISTA is than FISTA. 

The results in Tables \ref{table5} and \ref{table6} show that sFISTA is  more efficient than FISTA in all test problems. For example, when we focus on one row or one column of Table \ref{table5} or Table \ref{table6}, it is easy to see that all the values of $tratio$ are less than $1$. Moreover, $tratio$ in both tables are quite close to $0.2$, which indicates that the efficiency of sFISTA does not depend on the blurring type, the test images, the blurring levels or the noise types. These results further indicate that sFISTA is not only fast but also very robust for solving the image restoration problems considered in this paper.

\begin{table}[t]\footnotesize
\caption{Values of $tratio$ for \texttt{PRblurdefocus} with different blurring levels and types of noise}
\label{table5}
\begin{center}
\begin{tabular}{@{\extracolsep{\fill}}ccccccccccc}
\hline
\multirow{3}{*}{tratio} &\multirow{3}{*}{case1} &\multirow{3}{*}{case2} &\multirow{3}{*}{case3} &\multirow{3}{*}{case4} &\multirow{3}{*}{case5}&\multirow{3}{*}{case6}&\multirow{3}{*}{case7} &\multirow{3}{*}{case8} &\multirow{3}{*}{case9}\\
&&&&&&&&&\\
\multirow{3}{*}{$\mathrm{Blurlevel}$} &\multirow{3}{*}{`mild'} &\multirow{3}{*}{`mild'} &\multirow{3}{*}{`mild'} &\multirow{3}{*}{`medium'} &\multirow{3}{*}{`medium'}&\multirow{3}{*}{`medium'}&\multirow{3}{*}{`severe'} &\multirow{3}{*}{`severe'} &\multirow{3}{*}{`severe'}\\&&&&&&&&&\\
\multirow{3}{*}{Noise type} &\multirow{3}{*}{`gauss'} &\multirow{3}{*}{`laplace'} &\multirow{3}{*}{`multi'} &\multirow{3}{*}{`gauss'} &\multirow{3}{*}{`laplace'}&\multirow{3}{*}{`multi'}&\multirow{3}{*}{`gauss'} &\multirow{3}{*}{`laplace'} &\multirow{3}{*}{`multi'}\\
&&&&&&&&&\\
&&&&&&&&&\\
\hline
\multirow{3}{*}{`pattern1'} &\multirow{3}{*}{$0.1944$} &\multirow{3}{*}{$0.2016$} &\multirow{3}{*}{$0.1995$} &\multirow{3}{*}{$0.2371$} &\multirow{3}{*}{$0.2019$}&\multirow{3}{*}{$0.2052$}&\multirow{3}{*}{$0.2064$} &\multirow{3}{*}{$0.2027$} &\multirow{3}{*}{$0.1947$}\\
&&&&&&&&&\\
\multirow{3}{*}{`pattern2'} &\multirow{3}{*}{$0.1998$} &\multirow{3}{*}{$0.1971$} &\multirow{3}{*}{$0.1932$} &\multirow{3}{*}{$0.1982$} &\multirow{3}{*}{$0.2032$}&\multirow{3}{*}{$0.2007$}&\multirow{3}{*}{$0.2076$} &\multirow{3}{*}{$0.1979$} &\multirow{3}{*}{$0.1994$}\\
&&&&&&&&&\\
\multirow{3}{*}{`ppower'} &\multirow{3}{*}{$0.2021$} &\multirow{3}{*}{$0.2094$} &\multirow{3}{*}{$0.2373$} &\multirow{3}{*}{$0.1950$} &\multirow{3}{*}{$0.2009$}&\multirow{3}{*}{$0.1943$}&\multirow{3}{*}{$0.1994$} &\multirow{3}{*}{$0.1957$} &\multirow{3}{*}{$0.2022$}\\
&&&&&&&&&\\
\multirow{3}{*}{`smooth'} &\multirow{3}{*}{$0.2014$} &\multirow{3}{*}{$0.1983$} &\multirow{3}{*}{$0.1995$} &\multirow{3}{*}{$0.2068$} &\multirow{3}{*}{$0.1970$}&\multirow{3}{*}{$0.1999$}&\multirow{3}{*}{$0.1971$}&\multirow{3}{*}{$0.1965$} &\multirow{3}{*}{$0.2014$} \\
&&&&&&&&&\\
\multirow{3}{*}{`dot2'} &\multirow{3}{*}{$0.1968$} &\multirow{3}{*}{$0.1993$} &\multirow{3}{*}{$0.2171$} &\multirow{3}{*}{$0.2011$} &\multirow{3}{*}{$0.2024$}&\multirow{3}{*}{$0.1986$}&\multirow{3}{*}{$0.2028$} &\multirow{3}{*}{$0.2028$} &\multirow{3}{*}{$0.2074$}\\
&&&&&&&&&\\
\multirow{3}{*}{`dotk'} &\multirow{3}{*}{$0.2041$} &\multirow{3}{*}{$0.1995$} &\multirow{3}{*}{$0.1982$} &\multirow{3}{*}{$0.2021$} &\multirow{3}{*}{$0.2025$}&\multirow{3}{*}{$0.2143$}&\multirow{3}{*}{$0.1987$} &\multirow{3}{*}{$0.2016$} &\multirow{3}{*}{$0.2101$}\\
&&&&&&&&&\\
\multirow{3}{*}{`satellite'} &\multirow{3}{*}{$0.1987$} &\multirow{3}{*}{$0.1986$} &\multirow{3}{*}{$0.1989$} &\multirow{3}{*}{$0.2164$} &\multirow{3}{*}{$0.2088$}&\multirow{3}{*}{$0.2249$}&\multirow{3}{*}{$0.2046$} &\multirow{3}{*}{$0.1996$} &\multirow{3}{*}{$0.2049$}\\
&&&&&&&&&\\
\multirow{3}{*}{`hst'} &\multirow{3}{*}{$0.2006$} &\multirow{3}{*}{$0.1992$} &\multirow{3}{*}{$0.2027$} &\multirow{3}{*}{$0.1968$} &\multirow{3}{*}{$0.2016$}&\multirow{3}{*}{$0.2064$}&\multirow{3}{*}{$0.2165$} &\multirow{3}{*}{$0.2014$} &\multirow{3}{*}{$0.2034$}\\
&&&&&&&&&\\
&&&&&&&&&\\
\hline
\end{tabular}
\end{center}
\end{table}

\begin{table}[t]\footnotesize
\caption{Values of $tratio$ for \texttt{PRblurshake} with different blurring levels and types of noise}
\label{table6}
\begin{center}
\begin{tabular}{@{\extracolsep{\fill}}ccccccccccc}
\hline
\multirow{3}{*}{tratio} &\multirow{3}{*}{case1} &\multirow{3}{*}{case2} &\multirow{3}{*}{case3} &\multirow{3}{*}{case4} &\multirow{3}{*}{case5}&\multirow{3}{*}{case6}&\multirow{3}{*}{case7} &\multirow{3}{*}{case8} &\multirow{3}{*}{case9}\\
&&&&&&&&&\\
\multirow{3}{*}{$\mathrm{Blurlevel}$} &\multirow{3}{*}{`mild'} &\multirow{3}{*}{`mild'} &\multirow{3}{*}{`mild'} &\multirow{3}{*}{`medium'} &\multirow{3}{*}{`medium'}&\multirow{3}{*}{`medium'}&\multirow{3}{*}{`severe'} &\multirow{3}{*}{`severe'} &\multirow{3}{*}{`severe'}\\
&&&&&&&&&\\
\multirow{3}{*}{Noise type} &\multirow{3}{*}{`gauss'} &\multirow{3}{*}{`laplace'} &\multirow{3}{*}{`multi'} &\multirow{3}{*}{`gauss'} &\multirow{3}{*}{`laplace'}&\multirow{3}{*}{`multi'}&\multirow{3}{*}{`gauss'} &\multirow{3}{*}{`laplace'} &\multirow{3}{*}{`multi'}\\
&&&&&&&&&\\
&&&&&&&&&\\
\hline
\multirow{3}{*}{`pattern1'} &\multirow{3}{*}{$0.1984$} &\multirow{3}{*}{$0.2024$} &\multirow{3}{*}{$0.2172$} &\multirow{3}{*}{$0.1961$} &\multirow{3}{*}{$0.1984$}&\multirow{3}{*}{$0.1958$}&\multirow{3}{*}{$0.1973$} &\multirow{3}{*}{$0.2027$} &\multirow{3}{*}{$0.2138$}\\
&&&&&&&&&\\
\multirow{3}{*}{`pattern2'} &\multirow{3}{*}{$0.1995$} &\multirow{3}{*}{$0.2054$} &\multirow{3}{*}{$0.2167$} &\multirow{3}{*}{$0.1986$} &\multirow{3}{*}{$0.1974$}&\multirow{3}{*}{$0.1992$}&\multirow{3}{*}{$0.2010$} &\multirow{3}{*}{$0.1996$} &\multirow{3}{*}{$0.2219$}\\
&&&&&&&&&\\
\multirow{3}{*}{`ppower'} &\multirow{3}{*}{$0.2001$} &\multirow{3}{*}{$0.2028$} &\multirow{3}{*}{$0.2155$} &\multirow{3}{*}{$0.2081$} &\multirow{3}{*}{$0.2027$}&\multirow{3}{*}{$0.2195$}&\multirow{3}{*}{$0.2083$} &\multirow{3}{*}{$0.2113$} &\multirow{3}{*}{$0.2348$}\\
&&&&&&&&&\\
\multirow{3}{*}{`smooth'} &\multirow{3}{*}{$0.2025$} &\multirow{3}{*}{$0.2015$} &\multirow{3}{*}{$0.2022$} &\multirow{3}{*}{$0.2061$} &\multirow{3}{*}{$0.2083$}&\multirow{3}{*}{$0.2223$}&\multirow{3}{*}{$0.2168$}&\multirow{3}{*}{$0.2022$} &\multirow{3}{*}{$0.2169$} \\
&&&&&&&&&\\
\multirow{3}{*}{`dot2'} &\multirow{3}{*}{$0.1953$} &\multirow{3}{*}{$0.2052$} &\multirow{3}{*}{$0.2107$} &\multirow{3}{*}{$0.2141$} &\multirow{3}{*}{$0.1996$}&\multirow{3}{*}{$0.2034$}&\multirow{3}{*}{$0.2244$} &\multirow{3}{*}{$0.2350$} &\multirow{3}{*}{$0.1972$}\\
&&&&&&&&&\\
\multirow{3}{*}{`dotk'} &\multirow{3}{*}{$0.2114$} &\multirow{3}{*}{$0.2083$} &\multirow{3}{*}{$0.2148$} &\multirow{3}{*}{$0.2125$} &\multirow{3}{*}{$0.2051$}&\multirow{3}{*}{$0.2189$}&\multirow{3}{*}{$0.2049$} &\multirow{3}{*}{$0.2141$} &\multirow{3}{*}{$0.2081$}\\
&&&&&&&&&\\
\multirow{3}{*}{`satellite'} &\multirow{3}{*}{$0.2546$} &\multirow{3}{*}{$0.2129$} &\multirow{3}{*}{$0.2167$} &\multirow{3}{*}{$0.1965$} &\multirow{3}{*}{$0.2550$}&\multirow{3}{*}{$0.2123$}&\multirow{3}{*}{$0.2118$} &\multirow{3}{*}{$0.2210$} &\multirow{3}{*}{$0.2206$}\\
&&&&&&&&&\\
\multirow{3}{*}{`hst'} &\multirow{3}{*}{$0.1991$} &\multirow{3}{*}{$0.2047$} &\multirow{3}{*}{$0.2233$} &\multirow{3}{*}{$0.2071$} &\multirow{3}{*}{$0.1984$}&\multirow{3}{*}{$0.2577$}&\multirow{3}{*}{$0.2026$} &\multirow{3}{*}{$0.2025$} &\multirow{3}{*}{$0.2160$}\\
&&&&&&&&&\\
&&&&&&&&&\\
\hline
\end{tabular}
\end{center}
\end{table}
\end{exmp}

\section{Conclusion}
\label{sec:Conclusion}
In this paper, we propose the \emph{structured fast iterative shrinkage-thresholding algorithm} (sFISTA) for solving large scale ill-posed linear inverse problems arising from image restoration. By exploiting both the Kronecker product structure of the coefficient matrix and the pattern structure of the matrices in the Kronecker product approximation, sFISTA can significantly accelerate the computation compared to FISTA.  A theoretical error analysis has been conducted to show that sFISTA can reach the same level of computational accuracy as FISTA under certain conditions. Finally, the efficiency of sFISTA is demonstrated with both a theoretical computational complexity analysis and various numerical examples coming from different applications.  

The proposed sFISTA framework provide the possibility of developing new solvers for other imaging deburring problems. For example, it is possible to adapt sFISTA to solve nonsmooth optimization problems with sparsity constraints such as those using $l_1$-based regularization. We also plan to exploit preconditioning techniques to further reduce the iteration number and iteration time in our future work. 

\clearpage

\bibliographystyle{siamplain}
\bibliography{references}

\end{document}